\documentclass[12pt,twoside]{amsart}
\usepackage{amssymb}

\nonstopmode

\numberwithin{equation}{section} 

\font\tengothic=eufm10 scaled\magstep 1
\font\sevengothic=eufm7 scaled\magstep 1
\newfam\gothicfam
      \textfont\gothicfam=\tengothic
      \scriptfont\gothicfam=\sevengothic

\newtheorem{theorem}{Theorem}[section]
\newtheorem{lemma}[theorem]{Lemma}
\newtheorem{proposition}[theorem]{Proposition}

\theoremstyle{definition}

\newtheorem{definition}[theorem]{Definition} % \theoremstyle{remark}
\newtheorem{remark}[theorem]{Remark}
\newtheorem{example}[theorem]{Example}

\newtheorem{notation and remark}[theorem]{Notation and Remark}
\newtheorem{definition and remark}[theorem]{Definiton and Remark}
\newtheorem{remark and definition}[theorem]{Remark and Definition}

\newcommand\rank{\operatorname{rank}}

\newcommand{\al}{\alpha}
\newcommand{\be}{\beta}

\newcommand{\cH}{{\mathcal H}}
\newcommand{\cW}{{\mathcal W}}

\newcommand{\kindalong}{\ \hbox{$\hbox to .35in{\rightarrowfill}$} \ }
\newcommand{\reallylong}{\hbox{$\hbox to .5in{\rightarrowfill}$}  }

\def\ra{\rightarrow}

\def\End{{\rm End}}
\def\In{{\rm In}}

\def\ba{{\bf a}}
\def\bff{{\bf f}}
\def\Min{{\rm Min}}
\def\ra{\longrightarrow}

\def\deg{{\rm deg }}
\def\End{{\rm End}}

\def\rank{{\rm rank}}

\def\Ker{{\rm Ker}}
\def\ker{{\rm ker}}

\def\ol{\overline}
\def\ti{\tilde}
\def\noi{\noindent}

\def\cV{\mathcal{V}}
\def\cU{\mathcal{U}}
\def\cW{\mathcal{W}}

\def\Gr{{\rm Gr}}

\def\al{\alpha}
\def\be{\beta}

\def\In{{\rm In}}

\def\Min{{\rm Min}}
\def\ra{\longrightarrow}

\def\deg{{\rm deg }}
\def\End{{\rm End}}

\def\rank{{\rm rank}}

\def\Ker{{\rm Ker}}
\def\Max{{\rm Max}}
\def\bff{{\rm f}}

\def\Gr{{\rm Gr}}
\def\bff{{\bf f}}

\newfont{\bg}{cmr10 scaled\magstep4}
\newcommand{\bigzerol}{\smash{\hbox{\bg 0}}}
\newcommand{\bigzerou}{\smash{\lower1.7ex\hbox{\bg 0}}}
\newcommand{\bigast}{\smash{\lower3ex\hbox{\bg *}}}

\begin{document}

%%%%%%%%%%%%%%%%%%%%%%%%%%%%%%%%%%%%%%%%%%%%%%%%%%%%%%%%%%%%%%%%%%%%%%%%%%%

\begin{center}
{\Large \bf 
The strong Lefschetz property for Artinian algebras \\[1ex]
with non-standard grading
}
\end{center}
\medskip\medskip\medskip

\begin{center}
{\sc Tadahito Harima} \\[1ex]
Department of Mathematics, Hokkaido University of Education, \\
Kushiro  085-8580, Japan \\
E-mail: harima@kus.hokkyodai.ac.jp 
\medskip\medskip

{\sc Junzo Watanabe} \\[1ex]
Department of Mathematics, Tokai University, \\
Hiratsuka 259-1292, Japan \\
E-mail: junzowat@keyaki.cc.u-tokai.ac.jp
\end{center}
\medskip\medskip\medskip

%\maketitle

%%%%%%%%%%%%%%%%%%%%%%%%%%%%%%%%%%%%%%%%%%%%%%%%%%%%%%%%%%%%%%%%%%%%%%%%%
%%%%%%%%%%%%%%%%%%%%%%%%%%%%%%%%%%%%%%%%%%%%%%%%%%%%%%%%%%%%%%%%%%%%%%%%%

\begin{quote}
{\footnotesize 
{\sc Abstract.}\ 
Let $A=\oplus_{i=0}^cA_i$ be a graded Artinian $K$-algebra, 
where $A_c\neq(0)$ and ${\rm char} K=0$. 
({\em The grading may not necessarily be standard.}) 
Then $A$ has the strong Lefschetz property 
if there exists an element $g \in A_1$  such that 
the multiplication $\times g^{c-2i}:A_i \ra A_{c-i}$ is bijective 
for every $i=0,1,\ldots,[c/2]$. 
The main results obtained in this paper are as follows: 
\begin{itemize}
\item[1.] 
$A$ has the strong Lefschetz property 
if and only if there is a linear form $z\in A_1$ such that 
$\Gr_{(z)}(A)$ has the strong Lefschetz property. 
\item[2.] 
If $A$ is Gorenstein, then $A$ has the strong Lefschetz property if and only if there is a linear 
form $z \in A$ such that  
all central simple modules of $(A, z)$ have the strong Lefschetz property. 
\item[3.] 
A finite free extension of an Artinian $K$-algebra 
with the strong Lefschetz property 
has  the strong Lefschetz property if the fiber does.   
\item[4.] 
The complete intersection defined by power sums of consecutive degrees 
has the strong Lefschetz property. 
\end{itemize}
}
\end{quote}

%\medskip\medskip
%\noindent
%{\sc Keywords}: 
%Lefschetz property, Artinian algebra, complete intersection, Gorenstein algebra
%central simple module, 
%}
%\end{quote}

%%%%%%%%%%%%%%%%%%%%%%%%%%%%%%%%%%%%%%%%%%%%%%%%%%%%%%%%%%%%%%%%%%%%%%%%%
%%%%%%%%%%%%%%%%%%%%%%%%%%%%%%%%%%%%%%%%%%%%%%%%%%%%%%%%%%%%%%%%%%%%%%%%%

%%%%%%%%%%%%%%%%%%%%%%%%%%%%%%%%%%%%%%%%%%
%%%%%%%%%%%%%%%%%%%%%%%%%%%%%%%%%%%%%%%%%%
%%%%%%%%%%%%%%%%%%%%%%%%%%%%%%%%%%%%%%%%%%
%%%%%%%%%%%%%%%%%%%%%%%%%%%%%%%%%%%%%%%%%%

%Main theorems are Theorem \ref{main_thm_one}  Theorem \ref{main-th11} and 
%Theorem \ref{main-th2} and Theorem \ref{main-th3}.  

%\input{introduction.tex}
%%%%%%%%%%%%%%%%%%%%%%%%%%%%%%%%%%%%%%%%%%%%%%%%%%%%%%%%%%%%%%
%
%
%%%%      End of Abstract --- start of Introduction     %%%%%%
%
%
%%%%%%%%%%%%%%%%%%%%%%%%%%%%%%%%%%%%%%%%%%%%%%%%%%%%%%%%%%%%%%

%Main theorems are  
%Theorem \ref{main-th11} and 
%Theorem \ref{main-th2} and 
%Theorem \ref{main-th3}.  
%  and \begin{theorem} \label{proposition3-10}   %thm3.9

\section{Introduction}

By now a considerable amount of results have been obtained for the strong and weak Lefschetz properties 
for Artinian graded $K$-algebras over a field $K$ of characteristic zero
(e.g. \cite{tHjMuNjW01}, \cite{aW03}). 
So far all algebras considered for the strong Lefschetz property have the standard grading.  
In this paper
  we consider the Lefschetz properties  for algebras with grading which may {\em not necessarily be standard}.

The definition of the strong or weak  Lefschetz property 
 can be made verbatim, provided that the algebra has a linear form,  as with the case of the standard grading.   
With this definition,  one  notices  that some basic facts  of the  
Lefschetz property  in the case of standard grading easily fail to generalize.  
This is because 
a non-standard grading can  force tough restrictions on the behavior of a linear element as a  multiplication operator on the algebra. 
By the same reason, however, 
there are situations where the strong Lefschetz property is easier to prove 
for algebras with a non-standard grading. 

To explain in more detail, let $R=K[x_1, \ldots, x_n]$ be the polynomial ring, and let 
$$I=(p_d, p_{d+1}, \ldots,  p_{d+n-1})$$ 
be the ideal of $R$, where $p_d$ is the power sum symmetric function of degree $d$, i.e., 
\[
p_d=x_1^d + \cdots + x_n^d. 
\]  
Since the generators of $I$ are symmetric forms, they are contained in the ring of invariants $S \subset R$ under the action 
of the symmetric group.  Notice  that  $S$, being generated by  the elementary symmetric polynomials, 
has no longer the standard grading as a graded  subring of $R$.  
Put $A=R/I$ and $B=S/I \cap S$.  
Since   $e_1=x_1 + \cdots + x_n$ is a unique linear form of $S$,  
it is the only candidate for a strong Lefschetz element for $S/I \cap S$.  
In this particular case it happens that it is in fact a strong Lefschetz element. 
It is easier to prove the strong Lefschetz property for 
$S/I \cap S$ rather than for $R/I$.   
Once it is known that $S/I \cap S$ has the strong Lefschetz property, it determines the 
central simple modules of $(R/I, e_1)$, and it reduces the problem to a lower dimensional case. 
Our results of this paper enables us to conclude that $R/I$ indeed has the strong Lefschetz property.  

Generally speaking, it may well occur that the algebra 
$S/I \cap S$ does not have the strong Lefschetz property even if $I$ can be generated by a regular sequence of invariant forms. 
For example if the ideal contains a power of $e_1$  as a member of a minimal generating set, then  obviously it cannot have 
the strong Lefschetz property.  Nonetheless it reduces the problem to prove the strong Lefschetz property  to a lower dimensional case at least.  

Much the same results  of \cite{tHjW03} and  \cite{tHjW06}  are  valid for algebras with any grading as long as 
the algebra  possesses  a linear form. 
This is not surprising as the strong or  weak Lefschetz property is a property of a single linear form as an operator for a 
graded vector space.     

It seems natural to define  the strong Lefschetz property  for  a finite  graded vector space (rather than for an Artinian $K$-algebra)  
$$V=\oplus _{i=a}  ^{b} V_i.$$   For $V$ 
we call $g  \in  \End(V)$  a strong Lefschetz element if it is a grade-preserving map of degree one  and if the restricted map  
$$g^{b-a-2i}| _{V_{a+i}} : V_{a+i} \rightarrow V_{b-i}$$ 
is bijective for all $i$ such that $0 \leq  2i \leq b-a$.  

If $A$ is  a graded algebra  and if $V$ is a graded module over $A$, then we have the regular representation 
$$\times : A \ra \End(V).$$
Namely,  $\times z$ is the map $V  \ra V$ defined by $a \mapsto za$ for $a \in A$. 
Then we will say that $V$ has the strong Lefschetz property as an $A$-module if there 
exists a strong Lefschetz element  for $V$  in   the set $\times A \subset \End(V)$.  
In this case we say that  a linear form $z \in A$ is a strong Lefschetz element for $V$ 
if $\times z$ is.       
This way it is possible to extend the definition of  the strong Lefschetz property to graded modules over graded algebras whose grading  may 
not be  standard.   
In our previous paper \cite{tHjW06} we defined the central simple modules,  
$$U_1, \ldots, U_s, $$ 
 of a pair $(A, z)$ of an Artinian  algebra $A$  and a linear form $z \in A$,  in which case $A$ was assumed to have 
the  standard grading.  
We show that  much the same results can be extended to Artinian algebras with non-standard grading.  
Thus it is possible to establish the chain  of implications: 
\[
A  \mbox{ has SLP } 
\Leftarrow \Gr_{(z)}(A) \mbox{ has SLP }    
\Leftarrow  \oplus \tilde{U_i} \mbox{ has SLP }  
\Leftarrow  U_i \mbox{ has SLP for $\forall$ } i, 
\]
\noi
where $A$ is an Artinian Gorenstein algebra and $z$ is any  linear form in $A$.  
We have the  same  implications for an arbitrary Artinian algebra $A$ with 
certain additional  conditions for $\ti{U_i}$.       
For definition of $U_i$ and $\tilde{U}_i$ see Section~5. 

More important than the individual results is 
the fact that they can be used to prove the strong Lefschetz property for Artinian algebras with the 
{\em standard} grading.  An example is the complete intersection defined by power sums of consecutive degrees as mentioned earlier. 
This we discuss in section~7.  In section~8 we show more examples where the consideration of a non-standard grading works effectively.  

In section 3, we prove  that  $V_1 \otimes  V_2$ has the strong Lefschetz property if and only if  $V_1$ and $V_2$ have the  
the strong Lefschetz  property, where $V_i$ are  finite  graded vector spaces (Theorem~\ref{main_thm_one}). 
The ``if'' part  is fairly well-known for algebras  and the ``only if'' part  does not seem to have  been  written anywhere.   
One finds that the ``if''  part,  in the end, reduces  to the simplest case which asserts  that 
the algebra $$K[X, Y]/(X^r, Y^s)$$ has the strong Lefschetz property with the linear form $X+Y$ as a strong Lefschetz element.    
This itself is by no means trivial as the assertion is  essentially equivalent to what is known 
as the Gordan-Clebsch decomposition of the tensor product 
$$V(r) \otimes V(s)$$ of irreducible modules 
$V(r)$ and $V(s)$ over the Lie algebra  $sl_2$.  
One other method to prove it is to use the fact  which says  that the 
generic initial ideal  is Borel fixed  (\cite{tHjMuNjW01}).  
For completeness we show a new method to prove it using a theorem due  to  Ikeda in Appendix. 
The proof in itself is of considerable interest.  
For other results of this paper the discussion of section~3  is indispensable.         

In section~4 we show that the strong Lefschetz property  of $\Gr_{(z)}(A)$ implies that of $A$, where 
$z$ is any linear form (Theorem~\ref{main-th11}).    

In section~5 we generalize the notion of the central simple modules for non-standard grading 
and characterize the strong Lefschetz property in terms of it in two theorems. 

In section 6 we show that a finite free extension 
of an Artinian algebra with the strong Lefschetz property has the strong Lefschtez property if the fiber does. 
It is stated in Theorem~\ref{main-th3}.  
For the standard grading this was proved in \cite{tHjW03} with a correction made in \cite{tHjW04}.  
In this paper it is proved for any grading and the proof is substantially simplified.  
This is  another example which shows that  the ``central simple modules'' are useful.  

Section 2 is preliminaries, where we briefly review basic definitions and re-produce proofs for  some lemmas to be used in the sequel.  

Some related results can be found in Maeno~\cite{tM06} and Morita-Wachi-Watanabe~\cite{hMaWjW05}, where one finds a vast amount of examples of 
graded vector spaces with a strong Lefschetz element. 
This was really the starting point of the present paper. 

%%%%%%%%%%%%%%%%%%%%%%%%%%%%%%%%%%%%%%%%%%%%%%%%%%%%%%%%%%%%%%
%
%
%%%%      End of Introduction --- start of Section2     %%%%%%
%
%
%%%%%%%%%%%%%%%%%%%%%%%%%%%%%%%%%%%%%%%%%%%%%%%%%%%%%%%%%%%%%%

\section{Hilbert functions and the Lefschetz properties}

The Hilbert function of a graded vector space $V=\oplus_{i=a}^bV_i$ 
is the map $i \mapsto \dim V_i$. 
If $V$ has finite dimension, 
then its Hilbert series is the polynomial 
$$
h_V(q)=\sum_{i=a}^b(\dim V_i)q^i. 
$$
We define the {\em Sperner number} of $V$ by 
$$
{\rm Sperner}(V)={\rm Max}\{\dim V_a, \dim V_{a+1}, \ldots, \dim V_b\}
$$
and {\em CoSperner number}  by 
$$
{\rm CoSperner}(V)=\sum _{i=a}^{b-1} {\rm Min}\{\dim V_i, \dim V_{i+1}\}.
$$
The Hilbert function $h_V(q)$  of $V=\oplus_{i=a}^bV_i$ 
(where $V_a\neq (0)$ and $V_b\neq (0)$)  
is symmetric if 
$\dim V_{a+i}=\dim V_{b-i}$ for all $i=0,1,\ldots,[(b-a)/2]$. 
Then we call the half integer $(a+b)/2$ 
the reflecting degree of $h_V(q)$.     
The Hilbert function  of $V$ is unimodal 
if there exists an integer $m$ ($a\leq m\leq b$) such that
$$
\dim V_a\leq\dim V_{a+1}\leq\cdots\leq\dim V_m
\geq\dim V_{m+1}\geq\cdots\geq\dim V_b. 
$$

\begin{definition} \label{def2-1}   %def2.1
Let $A=\oplus_{i=0}^c A_i$ be a graded Artinian $K$-algebra. 
Suppose that $V=\oplus_{i=a}^b V_i$ 
is a finite graded $A$-module with $V_a\neq(0)$ and $V_b\neq(0)$. 
\begin{itemize}
\item[(i)]
$V$ has the {\em weak Lefschetz property} (WLP) as an $A$-module 
if there is a linear form $g\in A_1$ such that 
the multiplication $\times g: V_i\ra V_{i+1}$ is 
either injective or surjective for all $i=a,a+1,\ldots,b-1$. 
\item[(ii)]
$V$ has the {\em strong Lefschetz property} (SLP) as an $A$-module 
if there is a linear form $g\in A_1$ such that 
the multiplication $\times g^{b-a-2i}: V_{a+i}\ra V_{b-i}$ 
is bijective for all $i=0,1,\ldots,[(b-a)/2]$. 
\end{itemize} 
\end{definition}

\begin{remark}  %rem2.2 
\begin{itemize}
\item[(1)] It is easy to see that 
if $A$ has the WLP and if $A$ has the standard grading, then 
the Hilbert function of $A$ is unimodal.  
If the grading of $A$ is not standard, the WLP of $A$ does not imply the unimodality of 
the Hilbert function of $A$. 
Whatever the grading of   $A$ is,   the WLP of an  
$A$-module $V$  does not imply the unimodality of  the Hilbert function.  
\item[(2)] 
If $V$ has a unimodal Hilbert function, 
then ${\rm CoSperner}(V)=\dim V-{\rm Sperner}(V)$. 
\item[(3)] 
If $V$ has the SLP, 
then the Hilbert function of $V$ is symmetric and unimodal. 
\end{itemize}
\end{remark}

The following is a characterization of a weak Lefschetz element.

\begin{lemma}\label{lemma2-1}  %lemma2.3  
Let $A=\oplus_{i=0}^cA_i$ be a graded Artinian $K$-algebra 
and let $V=\oplus_{i=a}^b V_i$ be a finite graded $A$-module
with a unimodal Hilbert function, 
where $V_a\neq(0)$ and $V_b\neq(0)$. 
Then the following are equivalent. 
\begin{enumerate}
\item[$(1)$]
A linear form $\ell \in A_1$ is a weak Lefschetz element of $V$. 
\item[$(2)$]
$\dim V/\ell V = {\rm Sperner}(V)$. 
\item[$(3)$] 
$\rank (\times \ell)={\rm CoSperner}(V)$. 
\end{enumerate}
\end{lemma}

\begin{proof}

First note that for  any linear form $\ell  \in A_1$ we have 
$\dim V/\ell V \geq {\rm Sperner}(V)$ 
or equivalently,  
$\rank (\times \ell)\leq{\rm CoSperner}(V)$.
Now one sees that the equality holds if and only if 
$\ell : V_i \ra V_{i+1}$ is either surjective or injective for every $i$. 
(cf. Proposition 3.2 in \cite{jW87a}) 
\end{proof}

With the same notation as Definition~\ref{def2-1}, put 
$$
{\rm SP}_k(V)=\sum_{i=a}^b\Max\{\dim V_i - \dim V_{i-k}, \ 0 \} 
$$
for $1\leq k\leq b-a$, where $\dim_K V_j=0$ for all $j<a$. 
We call 
$$
{\bf SP}(V)=({\rm SP}_1(V), {\rm SP}_2(V),\ldots,{\rm SP}_{b-a}(V))
$$
the {\em Sperner vector} of $V$. 
Note that ${\rm SP}_1(V)$ is equal to the Sperner number of $V$ 
when $V$ has a unimodal Hilbert function. 
\medskip

Using the Sperner vector we can characterize a strong Lefschetz element  as follows.

\begin{lemma} \label{lemma2-2}  %lemma2.4  
Let $A=\oplus_{i=0}^c A_i$ be a graded Artinian $K$-algebra 
and let $V=\oplus_{i=a}^b V_i$ be a finite graded $A$-module 
with a symmetric unimodal Hilbert function, 
where $V_a\neq(0)$ and $V_b\neq(0)$. 
Then  the following conditions are equivalent.  
\begin{enumerate}
\item[$(1)$]
A linear form $g \in A_1$ is a strong Lefschetz element of $V$. 
\item[$(2)$]
$
\dim V/g^k V = {\rm SP}_k(V)
$ 
\ for all $k=1,2,\ldots,b-a$. 
\item[$(3)$]
$
\rank (\times g^k)=\dim V-{\rm SP}_k(V)
$   
\ for all $k=1,2,\ldots,b-a$. 
\end{enumerate}
\end{lemma}

\begin{proof}
Note that  we have  
$\dim_K V/f V \geq {\rm SP}_k(V)$
or equivalently, 
$\rank (\times f)\leq\dim V-{\rm SP}_k(V)$  
\ for any homogeneous form $f$ of degree $k$,  
with $1\leq k\leq b-a$.  
Now suppose that $f=g^k$ for a linear form $g$.  
Then  it is easy to see that  the equality holds if and only if 
$g^k :  V_i \ra V_{i+k}$ has the full rank for all $i$, and that   
it  is quivalent to claiming that $g$ is a strong Lefschetz element.   
(cf. Lemma 2.2 in \cite{tHjW06}) 
\end{proof}

For later use we would like to summarize the basic facts  used in the proofs of the 
above lemmas.

\begin{lemma}   \label{lemma2.5}   %lem2.5
Let $A$ be a graded Artinian $K$-algebra and let $V=\oplus_{i=a}^b V_i$ be a finite graded $A$-module 
with a symmetric unimodal Hilbert function, where $V_a\neq(0)$ and $V_b\neq(0)$.  
Then 
\begin{enumerate}
\item
$\dim V/\ell V \geq {\rm Sperner}(V)$ for all  $\ell  \in A_1$.  
\item
$\rank (\times \ell)\leq{\rm CoSperner}(V)$ for all $\ell  \in A_1$.  
\item 
$\dim_K V/f V \geq {\rm SP}_k(V)$ for all $f  \in A_k$, where $1\leq k\leq b-a$.  
\item  
$\rank (\times f)\leq\dim V-{\rm SP}_k(V)$  for any $f \in A_k$, where  $1\leq k\leq b-a$.  
\end{enumerate}
\end{lemma}

%%%%%%%%%%%%%%%%%%%%%%%%%%%%%%%%%%%%%%%%%%%%%%%%%%%%%%%%%%%%%%%%%%%%%%%
%
%
%%%%           End of Section 2 --- start of Section3            %%%%%%
%
%
%%%%%%%%%%%%%%%%%%%%%%%%%%%%%%%%%%%%%%%%%%%%%%%%%%%%%%%%%%%%%%%%%%%%%%%

\section{The tensor product of modules with the strong Lefschetz property}

In this section, we discuss the tensor product of modules with the strong Lefschetz property.  An important result is
Proposition~\ref{lemma_SLP_WLP} which often enables us to reduce an issue of the strong Lefschetz property to that of the 
weak Lefschetz property.  
\medskip 
\begin{notation and remark} \label{notation and remark3-1}  %remark3.1   
Let $A$ be a graded Artinian $K$-algebra, 
let $V$ be a finite graded $A$-module 
and let $z$ be a linear form of $A$. 
Since $V$ is a  finite graded  $A$-module,    
the linear map $\times z \in \End(V)$ is nilpotent. 
Hence the Jordan canonical form $J$ of $\times z$ is the matrix of the  following form: 
$$
J=
\left[
\begin{array}{ccccccc}
J(0,n_1) &     &        & \bigzerou         \\
    & J(0,n_2) &        &                   \\
    &     & \ddots &                        \\
\bigzerol   &      &        & J(0,n_r)      \\     
\end{array}
\right],  
$$
where $J(0,m)$ is the Jordan block of size $m\times m$ 
$$
J(0,m)=
\left[
\begin{array}{ccccccc}
0 &  1   &        & \bigzerou       \\
    & 0 &  \ddots      &            \\
    &     & \ddots &    1           \\
\bigzerol   &      &        & 0     \\     
\end{array}
\right].  
$$
Then we denote the {\em Jordan decomposition} of $\times z$ by writing 
$$
P(\times z)=n_1\oplus n_2\oplus \cdots \oplus n_r. 
$$
We note that 
$r=\dim_K V/zV$ and $\dim V=\sum n_i$. 
 Two decompositions 
$n_1\oplus n_2\oplus \cdots \oplus n_r$ 
and $n_1'\oplus n_2'\oplus \cdots \oplus n_{r'}'$ 
are regarded  as the same 
if they are the same as multisets. 
\end{notation and remark}

\begin{lemma}\label{lemma3}   %lemma3.2 
Let $A$ and $B$ be graded Artinian $K$-algebras.
Let $V$ be a finite graded $A$-module and $g \in A$ a linear form,  and similarly 
$W$ and $h$ for $B$. 
Assume that 
\begin{itemize}
\item[(i)] 
the Jordan decomposition of $\times g \in \End(V)$ is 
the same as that of $\times h \in \End(W)$, 
\item[(ii)] 
there is an integer $m$ such that 
$h_V(q)=q^mh_W(q)$. 
\end{itemize}
Then $g$ is a weak (resp. strong) Lefschetz element of $V$ 
if and only if 
$h$ is a weak (resp. strong) Lefschetz element of $W$. 
\end{lemma}

\begin{proof}
Since $h_V(q)=q^mh_W(q)$, 
the Sperner vectors of $V$ and $W$ are the same. 
Hence, using the assumption (i),  
the  assertion follows from Lemma~\ref{lemma2-1} 
and Lemma~\ref{lemma2-2}. 
\end{proof}

Following is another characterization of a strong Lefschetz element. 

\begin{lemma} \label{lemma3-3}   %lemma3.3
Let $A$ be a graded Artinian $K$-algebra 
and let $V=\oplus_{i=a}^bV_i$ be a finite graded $A$-module, 
where $V_a\neq(0)$ and $V_b\neq(0)$. 
Put $r={\rm Sperner}(V)$.   Let $g$ be a linear form of $A$. 
Then the  following conditions are equivalent. 
\begin{itemize}
\item[(1)]
$g$ is a strong Lefschetz element for $V$. 
\item[(2)]
There are  graded vector subspaces $\cV_1,\cV_2,\ldots,\cV_r$ of $V$ 
with  $V=\oplus_{i=1}^r\cV_i$ which satisfy the following conditions for each $i=1,2,\ldots,r$. 
\begin{itemize}
\item[(i)]
$g\cV_i\subset\cV_i$.  
\item[(ii)]
The Jordan canonical form of $\times g \in \End(\cV_i)$ 
is a single Jordan block.   
\item[(iii)]
The reflecting degree of $h_{\cV_i}(q)$ 
is equal to that of $h_V(q)$. 
\end{itemize}
\end{itemize}
In this case %$g$ is a strong Lefschetz element of $V$ and 
$P(\times g)=
\dim\cV_1\oplus\dim\cV_2\oplus\cdots\oplus\dim\cV_r$. 
\end{lemma}

\begin{proof}
(1) $\Rightarrow$ (2) 
Assume $g$ is a strong Lefschetz element for $V$. 
Let $h_V(q)=\sum_{i=a}^bh_iq^i$ be the Hilbert function of $V$. 
A basis for the Jordan decomposition of 
$\times g \in \End(V)$ is obtained as follows. 
Let $v_1,v_2,\ldots,v_{h_a}$ be a basis of $V_a$. 
By the SLP of $V$, 
the elements 
$$
\{g^kv_j \mid 1\leq j\leq h_a, 0\leq k\leq b-a\}, 
$$
being linearly independent, 
will be a part of the basis. 
Next let 
$\{v_{h_a+1}',v_{h_a+2}',\ldots,v_{h_{a+1}}'\}$ 
be a basis of ${\rm Ker} [V_{b-1} \stackrel{\times g}{\ra} V_b]$ 
(if it exists). 
By the SLP of $V$, 
there exist elements 
$\{v_{h_a+1},v_{h_a+2},\ldots,v_{h_{a+1}}\}$ of $V_{a+1}$ 
such that 
$g^{b-a-2}v_j=v_j'$ for all $j=h_a+1,h_a+2,\ldots,h_{a+1}$. 
Then the elements 
$$
\{g^kv_j \mid h_a+1\leq j\leq h_{a+1}, 0\leq k\leq b-a-2\}, 
$$
none of these being dependent of the previously chosen basis elements, 
will be another part of the basis. 
We repeat the same  to expand basis elements. 
We may carry over this process to decompose $\times g$ into Jordan blocks. 
Since 
$$
h_a+(h_{a+1}-h_a)+(h_{a+2}-h_{a+1})+\cdots+
(h_{[(b-a)/2]}-h_{[(b-a)/2]-1})=r,   
$$
one sees that there are $r$ Jordan blocks.  

Now let $\cV_i$ be the subspace  of $V$ spanned by   
$\{v_i,gv_i,\ldots,g^{b-a-2(d_i-a)}v_i\}$
for   $i=1,2,\ldots,r$, 
where $d_i=\deg(v_i)$. 
Then it is easy to verify 
%%%%%%%%%%that $g$ and $\cV_i$ $(1\leq i\leq r)$ satisfy 
the conditions stated in (2). 

(2) $\Rightarrow$ (1) By (iii) it suffices to prove this for each $\cV _i$, which is obvious.    
\end{proof}

\begin{proposition} \label{proposition3.4}   %prop3.4
Let $K$ be a field of characteristic zero. 
Let $A$ and $B$ be graded Artinian $K$-algebras, 
let $V$  be a finite graded $A$-module  
and  $W$ a finite graded $B$-module. 
If $V$ and $W$ have the SLP, 
then $V\otimes_KW$ also has the SLP as an $A\otimes_KB$-module. 
\end{proposition}

First we prove  a lemma. 

\begin{lemma} \label{lemma3-5}   %lemma3.5 
With the same notation as Proposition~\ref{proposition3.4}, 
let $g \in A$ be any linear form   and let $\cV$  
be a graded subspace of $V$  such that   
$g \cV \subset \cV$.  
Similarly let $h \in B$ be a linear form  and $\cW \subset  W$ a graded  subsupace 
such that $h \cW \subset \cW$.  
Assume that $\times g \in \End(\cV)$ is a single Jordan block  
and  the same for $\times h \in \End(\cW)$.  
Moreover let $m=\dim \cV$,  $n=\dim \cW$ and $s={\rm Sperner}(\cV \otimes\cW)$.   

Then there exist graded vector subspaces $\cU_1,\cU_2,\ldots,\cU_s$ 
of $\cV \otimes \cW$ such that $\cV \otimes\cW =\oplus_{i=1}^s\cU_i$, 
which satisfy the following conditions for  $i=1,2,\ldots,s$.  
\begin{itemize}
\item[(i)]
$(g \otimes 1 + 1\otimes h)\cU_i\subset\cU_i$.   
\item[(ii)]
The Jordan canonical matrix of $\times(g \otimes 1+1\otimes h) \in \End(\cU_i)$ 
is a single Jordan block. 
\item[(iii)]
The reflecting degree of $h_{\cU_i}(q)$ 
is equal to that of $h_{\cV \otimes\cW}(q)$. 
\end{itemize} 
\end{lemma}

\begin{proof}
By  Lemma~\ref{lemma3-3}, 
it is enough to show that 
$g \otimes 1+1\otimes h$ is 
a strong Lefschetz element for  $\cV \otimes \cW$. 
Let $d \mbox{ and } e$ be the initial degrees of $\cV$ and $\cW$ respectively.  
Then the Hilbert functions  of $\cV$ and $\cW$ are: 
$$
\left\{\begin{array}{l}
h_{\cV}(q)=q^{d}+q^{d+1}+\cdots+q^{d+m-1},  \\
h_{\cW}(q)=q^{e}+q^{e+1}+\cdots+q^{e+n-1}. 
\end{array}
\right.
$$
Let $K[x]$ be the polynomial ring in one variable.   
As a graded vector space we may choose an isomorphims 
$\cV \ra K[x]/(x^m)(-d) $  so that we have the  commutative diagram: 

$$
\begin{array}{ccc}
\cV  & \stackrel{\times g}{\longrightarrow} & \cV \\[1ex]
\downarrow   &    &\downarrow    \\[1ex]
(K[x]/(x^{m}))(-d) 
& \stackrel{\times x}{\longrightarrow} & (K[x]/(x^{m}))(-d). 
\end{array}
$$
Likewise we may choose an isomorphism $\cW \ra K[y]/(y^n)(-e)$  
to  get the  similar diagram for $\cW$.  

These diagrams give rise  the following commutative diagram: 
$$
\begin{array}{ccc}
\cV \otimes\cW & 
\stackrel{\times (g \otimes 1+1\otimes h)}{\longrightarrow} 
& \cV \otimes\cW  \\[1ex]
\downarrow   &    &\downarrow    \\[1ex]
(K[x, y]/(x^{m}, y^{n}))(-(d+e)) 
& \stackrel{\times (x + y)}{\longrightarrow} 
& (K[x, y]/(x^{m}, y^{n}))(-(d + e)),  
\end{array}
$$ 
where the vertical maps are  isomorphisms as graded vector spaces. 
Since the characteristic of $K$ is zero, 
it follows by 
%(Corollary 3.5 in \cite{jW87a})  
Proposition~\ref{ikeda} in Appendix 
that $x + y$ is a strong Lefschetz element for 
$K[x, y]/(x^{m}, y^{n})$. 
Hence $g \otimes 1+1\otimes h$ is  
a strong Lefschetz element for  $\cV \otimes\cW$ as well as 
$x +  y$ is a strong Lefschetz element for  
$(K[x, y]/(x^{m}, y^{n}))$.  
\end{proof}

%%%%%%%%%%%%%%%%%%%%%%%%%%%%%%%%%%%%%%%%%%%%%%%%%%%%%%%%%%%%%%%%
\begin{proof}[Proof of Proposition~\ref{proposition3.4}]
Let $V=\oplus_{i=1}^r\cV_i$ and $W=\oplus_{i=j}^u\cW_j$ 
be the direct sum decomposition of $V$ and $W$ 
constructed in Lemma~\ref{lemma3-3} 
with respect to the strong Lefschetz elements
$g \in A_1$ and $h \in B_1$ for  $V$ and $W$, 
where $r={\rm Sperner}(V)$ and $u={\rm Sperner}(W)$. 
Then  we see that 
$V\otimes W=\oplus_{i,j}(\cV_i\otimes\cW_j)$ 
and
$\cV_i\otimes\cW_j$ is closed by 
the multiplication 
$\times(g \otimes 1+1\otimes h):V\otimes W\ra V\otimes W$ 
for all $i$ and $j$.  
Noting  
$$
\left\{\begin{array}{ll}
\ & h_{V\otimes W}(q)=h_V(q)h_W(q)
=(\sum_{i=1}^rh_{\cV_i}(q))(\sum_{j=1}^uh_{\cW_j}(q)) \\[1ex] 
%{\rm and} 
& 
h_{\cV_i\otimes \cW_j}(q)=h_{\cV_i}(q)h_{\cW_j}(q), 
\end{array}
\right. 
$$
we see that 
the reflecting degree of $h_{\cV_i\otimes\cW_j}(q)$ 
is equal to that of $h_{V\otimes W}(q)$. 
By Lemma~\ref{lemma3-5},  for each pair $(i,j)$, 
there are subspaces $\cU ^{(ij)}_k$ satisfying the conditions 
(i), (ii) and (iii) of Lemma~\ref{lemma3-5}.  
We note that the Sperner number of $V \otimes W$ is equal to the sum of those of 
$\cV _i \otimes \cW _j$.  Hence, by Lemma~\ref{lemma3-3},  the linear form $g \otimes 1 + 1 \otimes h$ 
is a strong Lefschetz element for $V \otimes W$.
\end{proof}

%%%%%%%%%%%%%%%%%%%%%%%%%%%%%%%%%%%%%%%%%%%%%%%%%%%%%%%%%%%%%

\begin{lemma} \label{dual_hilbert}  %3.6
Let \[
u_1 \geq u_2 \geq \cdots \geq u_r 
\]
 ba a descending sequence  of positive integrers, all even or all odd. 
Let $a,b$ be any integers such that $b-a+1 = u_1$.  
Let  

\begin{equation}   \label{xxx} 
h(q)=\frac{1}{1-q}  \sum _{i=1}^r  
\left(
q^{(a+b+1-u_i)/2}- q^{(a+b+1+u_i)/2}
\right)
\end{equation}
Then $h(q)$ is a symmetric unimodal polynomial with reflecting degree $(a+b)/2$ and the  coefficient  
$h_{a+k}$ of $q^{a+k}$ in $h(q)$ for $0 \leq  k \leq (b-a)/2$  is given by 
\[
h_{a+k}=\#\{u_i | u_i \geq  b-a + 1 - 2k\}.  
\]
Conversely suppose that 
Let $h(q)=\sum_{i=a}^bh_iq^i$ 
is a symmetric  unimodal polynomial with positive 
 integers as coefficients and $h_ah_b \not =0$.    
Put $r=\Max\{h_a,h_{a+1},\ldots,h_b\}$. 
Then $h(q)$ can be written as above 
with positive integers $u_1, \cdots, u_r$.  
These integers are uniquely determined as a multiset by $h(q)$.  
In particular $\mbox{\rm Max} \{u_i \} = b-a+1$ and $(a+b)/2$ is the 
reflecting degree of $h(q)$. 
\end{lemma} 
%\vspace{1ex}
\begin{proof}
Notice that the $i$-th summand (after divided by $1-q$) on the right hand side of~(\ref{xxx})     
is 
$$p_i(q):=q^{d_i}+q^{d_i+1}+\cdots+q^{d_i+(u_i-1)}, $$ 
where $d_i=(a+b+1-u_i)/2$.  
In particular, $p_1(q)=q^a+q^{a+1}+\cdots+q^b$.  Details are left to the reader.  
\end{proof}  

%%%%%%%%%%%%%%%%%%%%%%%%%%%%%%%%%%%%%%%%%%%%%%%%%%%%%%%%%%%%%%%%%%%%%%%%%%%%%%%%%%%%
\begin{lemma}  \label{lemma3.7}  %3.7 %new 
Suppose that $V$ is a finite graded module with the Hilbert function $h_V(q)$ 
over a graded Artinian algebra $A$.  Let $u_1, \ldots, u_r$ be the positive integers  
satisfying {\rm (\ref{xxx})} 
in the above lemma.  Let $g \in A_1$ and let $P(\times g)=n_1 \oplus \cdots \oplus n_{r'}$.  
Then $g$ is a strong Lefschetz element for $V$ if and 
only if $\{u_i\}$ and $\{n_i\}$ are the same as multisets.  
\end{lemma}

\begin{proof}
Assume that $g$ is a strong Lefschetz element for $V$.  
Let  $$\cV_1, \ldots, \cV _r$$  be the subspaces given in Lemma~\ref{lemma3-3}.   
Then  $\cV_1, \ldots, \cV _r$ 
correspond to Jordan blocks for $\times g \in \End(V)$.  Thus $n_i= \dim \cV _i$, if they are arranged in 
decreasing order.  By construction of the bases for $\cV_i$ one 
sees that $\{n_i\}$ and $\{u_i\}$ coincide as multisets. 
The converse follows from Lemma~\ref{lemma2-2}. 
\end{proof} 

%%%%%%%%%%%%%%%%%%%%%%%%%%%%%%%%%%%%%%%%%%%%%%%%%%%%%%%%%%%%%%%%%%%%%%%%%%%%%%%%%%%%

\begin{lemma} \label{lemma3.8}   %lemma3.8
Let $A$ and $B$ be graded Artinian $K$-algebras, where ${\rm char} (K)=0$.  
Let $V$ be a  finite graded   $A$-module 
and $W$ a finite graded $B$-module. 
\begin{enumerate}
\item
Let $z\in A_1$ and $z'\in B_1$.   
If $P(\times z)=d_1\oplus d_2\oplus\cdots\oplus d_r$ 
and $P(\times z')=f_1\oplus f_2\oplus\cdots\oplus f_{r'}$, 
then
$$
P(\times(z\otimes 1+1\otimes z'))
=\oplus_{i,j}\oplus_{k=1}^{\Min\{d_i,f_j\}}(d_i+f_j+1-2k). 
$$ 
In particular 
$$
{\rm corank}\left( \times (z\otimes 1 + 1 \otimes z')\right) = \sum _{i,j}\Min\{d_i,f_j\}. 
$$
\item
Let $\{u_i\}$ be the integers for $h_V(q)$ given in {\rm Lemma~\ref{dual_hilbert}} and similarly 
$\{v_j\}$ for $W$. 
Then 
$$
\mbox{\rm Sperner}(V \otimes  W)=\sum _{i,j}\Min\{u_i,v_j\}. 
$$
\end{enumerate}
\end{lemma}
\vspace{1ex}

\begin{proof}
(1) This can be proved  in  the same way 
as Proposition~10 in \cite{tHjW03}. 
(2) Sperner($V\otimes W$) is equal to the corank of the linear map 
$\times(z\otimes 1+1\otimes z')$, provided that it is a strong Lefschetz element.  
This is the case if $P(\times z)= u_1 \oplus \cdots \oplus u_r$ and 
$P(\times z')=v_1 \oplus \cdots \oplus v_{r'}$, where $r$ and $ r'$ 
are the Sperner numbers of $V$ and $W$ respectively.   
\end{proof}

\begin{proposition} \label{lemma_SLP_WLP}   %lemma3.9 
Let $K$ be a field of characteristic zero, 
let $A$ be a graded Artinian $K$-algebra 
and let $V=\oplus_{i=a}^bV_i$ be a finite graded $A$-module 
with a symmetric unimodal Hilbert function, 
where $V_a\neq(0)$ and $V_b\neq(0)$ and $b-a \geq 1$. 
Let $t$ be a new variable. 
Then the following conditions are equivalent. 
\begin{enumerate}
\item[$(1)$]
$V$ has the SLP. 
\item[$(2)$]
$V\otimes_K K[t]/(t^\alpha)$ has the WLP as an $A\otimes_K K[t]$-module 
for all positive integers $\alpha$. 
\item[$(3)$]
$V\otimes_K K[t]/(t^\alpha)$ has the WLP as an $A\otimes_K K[t]$-module 
for all $\alpha=1,2,\ldots,b-a$. 
\end{enumerate}
\end{proposition}

\begin{proof}
(1) $\Rightarrow$ (2) follows from Proposition~\ref{proposition3.4}.  
(2) $\Rightarrow$ (3) is trivial. 
We show (3) $\Rightarrow$ (1).
Put $r={\rm Sperner}(V)$.  
Let $u_1,u_2,\ldots,u_r$ be the decreasing sequence of positive integers for $h_V(q)$ 
given in  Lemma~\ref{dual_hilbert}. 
By way of contradiction assume  that 
no linear form of $A$ is a strong Lefschetz element for $V$.  
Let  $z\in A_1$ be any  element.    Suppose that  $\times z \in \End(V)$   
decomposes as   $P(\times z)=n_1\oplus n_2\oplus\cdots\oplus n_{r'}$, 
where $n_1\geq n_2\geq\cdots\geq n_{r'}$. 
By Lemma~\ref{lemma3.7},  the two multisets  $\{u_i\}$ and $\{n_j\}$ are 
different  since $\times z$ is not a strong Lefschetz element.  

Let $j$ be the least integer for which  $u_j \not = n_j$.  
We claim that $u_j > n_j$.  If $j=1$ or $r=1$, then the claim is obvious.     
Assume that $n_1=u_1$.  
Let $Z$ be a variable.  We may regard  $V$  as a graded module over $K[Z]$ via the algebra map 
$K[Z] \ra A$ defined by $Z \mapsto z$.  
Since $n_1=u_1=b-a+1$, there is a Jordan block  $\cV  \subset V$ for $\times z \in \End(V)$ 
such that the Hilbert function  is 
$$
h_{\cV}(q)=(q^a+q^{a+1}+\cdots+q^b)  %=\sum_{i=a'}^{b'}h_i'q^i, 
=\frac{1}{1-q}\left(q^{(a+b+1-u_1)/2} - q^{(a+b+1+u_1)/2}\right).  
$$
Let $V' \subset V$  be a $K[Z]$-module such that $V=V' \oplus \cV$ and such that the Jordan 
decomposition of $\times z \in \End(V')$ is given by  
 $n_2 \oplus \cdots \oplus n_{r'}$.  
Since the Hilbert function of 
$V'$ is $h_V(q)- h_{\cV}(q)$, which is symmetric unimodal with the maximal part one less 
than that of $h_V(q)$,  proof of the  claim is now complete  by induction of $r$.   

Now choose the least $j$ such that $u_j > n_j$.  
Put $\alpha=n_j$. 
Since $u_1=b-a+1$ and $u_1\geq u_2\geq\cdots\geq u_j>n_j=\alpha$, 
we see $\alpha\leq b-a$. 
Consider the  $A \otimes K[t]$-module 
$$W=V\otimes K[t]/(t^\alpha).$$ 
By Lemma~\ref{lemma3.8} (1),  
the Jordan decomposition of 
$\times(z\otimes 1+1\otimes t) \in  \End(W)$ is  given by 
$$
P(\times(z\otimes 1+1\otimes t))
=\oplus_{i=1}^{r'}[\oplus_{k=1}^{\Min\{n_i,\alpha\}}(n_i+\alpha+1-2k)]. 
$$ 
This implies 
$$\dim W/(z\otimes 1+1\otimes t)W =  
\sum_{i=1}^{r'}\Min\{n_i,\alpha\} \\[1ex].$$

On the other hand,   
since $\sum_{i=1}^{r'} n_i=\sum_{i=1}^ru_i$, 
we have that $\sum_{i=j+1}^{r'} n_i>\sum_{i=j+1}^r u_i$. 
Hence 
$$
\begin{array}{rcl}
\dim W/(z\otimes 1+1\otimes u)W & = & 
\sum_{i=1}^{r'}\Min\{n_i,\alpha\} \\[1ex]
  & = & 
\underbrace{\alpha+\cdots+\alpha}_{j}+n_{j+1}+\cdots+n_{r'} \\[1ex]
  & > & 
\underbrace{\alpha+\cdots+\alpha}_{j}+u_{j+1}+\cdots+u_r \\[1ex]
  & \geq & 
\sum_{i=1}^r\Min\{u_i,\alpha\}.  \\[1ex]
\end{array}
$$
By Lemma~\ref{lemma3.8} (2),  
we have  
$$
\sum_{i=1}^r\Min\{u_i,\alpha\} 
={\rm Sperner} (W). 
$$
Thus we have shown that no linear form of $A \otimes K[t]/(t^{\al})$ is a weak Lefschetz 
element for $W$, which contradicts  the assumption (3).   
\end{proof}

\begin{theorem} \label{main_thm_one}   %thm3.10
Let $K$ be a field of characteristic zero. 
Let $A$ and $B$ be  graded Artinian $K$-algebras, 
and let $V$ be a finite graded  $A$-module 
and $W$ a  finite graded  $B$-module. 
Suppose that 
the Hilbert functions of $V$ and $W$ are 
symmetric and unimodal.  Then 
 $V\otimes_KW$ has the SLP as an $A\otimes_KB$-module, 
if and only if  $V$ and $W$ have the SLP. 
\end{theorem}

\begin{proof}
The ``if'' part was proved in Proposition~\ref{proposition3.4}.
We prove the ``only if'' part.  
By way of contradiction assume that $V$ does not have the SLP. 
By Proposition~\ref{lemma_SLP_WLP}, the module  
$V^\prime=V\otimes K[t]/(t^n)$ does not have the WLP 
for some integer $n>0$. 

Let $v_1,v_2,\ldots,v_s$ be the integers for $h_{W}(q)$ 
given in  Lemma~\ref{dual_hilbert}. 
We may assume that $v_s=1$. 
In fact, if $v_s>1$, 
we may replace $W$ by $W\otimes K[u]/(u^m)$ 
where $m$ is the number of maximal part  
in the Hilbert function of $W$. 

Let $u_1,u_2,\ldots,u_r$ be the integers for $h_{V^\prime}(q)$ 
given in  Lemma~\ref{dual_hilbert}. 
Let $z$ be any linear form of $A\otimes K[t]$ 
and $P(\times z)=d_1\oplus d_2\oplus\cdots\oplus d_\alpha$ 
the Jordan decomposition of $\times z \in \End(V^\prime)$. 
Here note that $r \leq \alpha$. 
Also note that if $V' = \oplus _{1=1}^{\al} \cV _i$ is a decomposition into  Jordan blocks, then 
the Hilbert function of  $\cV _i$  is of  the  form  
$$q^a + q^{a+1}+ \cdots + q^{b}$$
for some integers $a, b$ (depending on $i$).   
Bearing this in mind, we  see  that 
%$$
\begin{eqnarray*}
u_1+u_2+\cdots+u_k & = & 
\sum_{i=-\infty} ^{\infty} 
\Min\{\mbox{$k$, the coefficient of  $q^i$  in   $h_{V^\prime}(q)$}\} \\[1ex]
  & \geq & d_1+d_2+\cdots+d_k
\end{eqnarray*}
%$$ 
for all $k=1,2,\ldots,r$. 
Hence, noticing that  $\sum_{i=1}^r u_i=\sum_{i=1}^\alpha d_i$, 
we see 
$$
u_k+u_{k+1}+\cdots+u_\alpha
\leq d_k+d_{k+1}+\cdots+d_\alpha 
$$ 
for all $k=1,2,\ldots,\alpha$, 
where we have put $u_i=0$ for  $i>r$. 
Using this, we can easily verify  
\begin{equation}\label{EQ:3-1}
\sum_{i=1}^r\Min\{u_i,f\} \leq \sum_{i=1}^\alpha\Min\{d_i,f\} 
\end{equation}
for any integer $f>0$. 
Let $y$ be any linear form of $B$ 
and $P(\times y)=f_1\oplus f_2\oplus\cdots\oplus f_\beta$ 
the Jordan decomposition of $\times y \in \End(W)$. 
Similarly we get  
\begin{equation}\label{EQ:3-2}
\sum_{j=1}^s\Min\{d, v_j\} \leq \sum_{j=1}^\beta\Min\{d,f_j\} 
\end{equation}
for any integer $d>0$. 

Recall that $r$ is the Sperner number of $V'$ and $\al$ is 
equal to the corank of the linear map  $\times z : V' \ra V'$. 
So we have 
$r<\alpha$ by Lemma~\ref{lemma2-1}  % from Lemma~\ref{lemma3-7} (3). 
as $z$ is not a weak Lefschetz element of $V^\prime$.   
Using the inequalities (\ref{EQ:3-1}) and (\ref{EQ:3-2}) above, 
we obtain
$$
\begin{array}{rcl}
\sum_{i,j}\Min\{u_i, v_j\} & = & 
\{\sum_{j=1}^{s-1}(\sum_{i=1}^r\Min\{u_i, v_j\})\} 
+ \sum_{i=1}^r\Min\{u_i, 1\} \\[1ex]
 & < & 
\{\sum_{j=1}^{s-1}(\sum_{i=1}^\alpha\Min\{d_i, v_j\})\} 
+ \sum_{i=1}^\alpha\Min\{d_i, 1\} \\[1ex]
& = & \sum_{i,j}\Min\{d_i, v_j\} \\[1ex]
& \leq & \sum_{i,j}\Min\{d_i, f_j\}. 
\end{array}
$$
This shows that 
$z\otimes 1+1\otimes y$ is not 
a weak Lefschetz element of $V^\prime\otimes W$, since  
$$
{\rm Sperner\ }(V^\prime\otimes W)=\sum_{i,j}\Min\{u_i, v_j\}
$$
and since 
$$
\dim (V^\prime\otimes W)/
(z\otimes 1+1\otimes y)(V^\prime\otimes W)
=\sum_{i,j}\Min\{d_i, f_j\}.   
$$

This means that $V^\prime\otimes W$ does not  have the WLP. 
However, since $V\otimes W$ has the SLP,  
Proposition~\ref{lemma_SLP_WLP} implies  that 
$V^\prime\otimes W$ has the WLP. 
This is a contradiction. 
 \end{proof}

%%%%%%%%%%%%%%%%%%%%%%%%%%%%%%%%%%%%%%%%%%%%%%%%%%%%%%%%%%%%%%%%%%%%%%%%%%%%
%%%%%%%%%%%%%%%%%%%%%%%%%%%%%%%%%%%%%%%%%%%%%%%%%%%%%%%%%%%%%%%%%%%%%%%
%
%
%%%%           End of Section 3 --- start of Section4           %%%%%%
%
%
%%%%%%%%%%%%%%%%%%%%%%%%%%%%%%%%%%%%%%%%%%%%%%%%%%%%%%%%%%%%%%%%%%%%%%%

\section{The strong Lefschetz property for $\Gr_{(z)}(A)$}   %section4
We need some preparations for the proof of Theorem~\ref{main-th11}. 

\begin{notation and remark} \label{notation and remark 3-1}  %remark4.1
Let $A=\oplus_{i=0}^cA_i$ be a graded Artinian $K$-algebra. 
For any linear form $z \in A_1$, 
consider the associated graded ring 
$$
\Gr_{(z)}(A)=A/(z)\oplus (z)/(z^2) \oplus (z^2)/(z^3) \oplus \cdots \oplus (z^{p-1})/(z^p), 
$$
where $p$ is the least integer such that $z^p  =0$. 
For a non-zero element $a\in A$ 
there is $i$ such that $a\in(z^i)\setminus  (z^{i+1})$. 
In this case we write $a^\ast\in\Gr_{(z)}(A)$ 
for the natural image of $a$ in $(z^i)/(z^{i+1})$. 

As is well known $\Gr_{(z)}(A)$ is endowed with a commutative ring structure. 
The multiplication in $\Gr_{(z)}(A)$ is given by 
$$
(a+(z^{i+1}))(b+(z^{j+1}))=ab+(z^{i+j+1}), 
$$
where $a \in (z^i) \setminus (z^{i+1})$ and $b \in (z^j)\setminus (z^{j+1})$.    
Note that 
$$
(z^i)/(z^{i+1}) \cong 
z^iA_0\oplus(z^iA_1/z^{i+1}A_0)\oplus  (z^iA_2/z^{i+1}A_1)\oplus \cdots 
\oplus(z^iA_{c-i}/z^{i+1}A_{c-i-1})
$$
as graded vector spaces for all $i=0,1,\ldots,p-1$. 
Furthermore note that $\Gr_{(z)}(A)$ inherits a grading from $A$. 
More precisely, $\Gr_{(z)}(A)=\oplus_{i=0}^c[\Gr_{(z)}(A)]_i$, where 
$$
[\Gr_{(z)}(A)]_i \cong
(A_i/zA_{i-1})\oplus(zA_{i-1}/z^2A_{i-2})\oplus\cdots
\oplus(z^{i-1}A_{1}/z^{i}A_{0})\oplus z^iA_0
$$
as graded vector spaces for all $i=0,1,\ldots,c$. 
Hence, 
$\Gr_{(z)}(A)$ and $A$ have the same Hilbert function. 
\end{notation and remark}

\begin{notation and remark} \label{notation and remark 3-2}   %remark4.2
Let $S=K[Y,Z]$ be the polynomial ring in two variables over an infinite field $K$. 
Here we regard S as a standard graded $K$-algebra 
with $\deg(Y)=\deg(Z)=1$.  
Let $V$ be a finite graded $S$-module. 
Write $V\cong F/N$, 
where 
$$
F=\oplus_{j=1}^s [ K[Y,Z](-d_j) ]
$$ 
is a free graded $S$-module of rank $s$ 
and $N$ a graded submodule of $F$. 
An element $\bff \in F$ can be written uniquely as 
$\bff=\ba_0+\ba_{1}Z+\cdots +\ba_dZ^{d}$ with $\ba_i \in \oplus^s K[Y]$.
Denote by $\In'(\bff)$ the term $\ba_j Z^j$ 
for the minimal $j$ such that $\ba _j \not  = {\bf 0}$. 
Furthermore we define 
${\rm In}'(N)$ to be the graded submodule of $F$ 
generated by the set $\{ {\rm In}' (\bff) \}$, 
where $\bff$ runs over homogeneous forms of $N$. 
Put  
$$
\Gr_{(Z)}(V)=V/ZV \oplus ZV/Z^2V \oplus Z^2V/Z^3V \oplus \cdots.  
$$
Then, we have that $\Gr_{(Z)}(V) \cong F/\In'(N)$ 
as finite graded $S$-modules. 
Suppose that 
$$
F_1 \stackrel{\phi} {\ra} F_0 \ra V \ra 0
$$ 
is a finite presentation of $V$.   
Let $\Delta_{\rm \tiny MAX}(V)$ be the ideal of $S$ 
generated by the maximal minors of $\phi$. 
As a Fitting ideal it does not depend on the choice of finite presentation. 
\end{notation and remark}

The following is a key to the proof of Theorem \ref{main-th11}.

%---------------------------------------
\begin{proposition}[\cite{tHjW06}  Proposition 3.3] \label{prop}     %prop4.3
%---------------------------------------
With the same notation as above, 
let $g$ be a general linear form of $S$. 
Assume that $K$ is an infinite field. 
Then 
$$
\dim _K V/gV \leq \dim _K \Gr_{(Z)}(V)/g\Gr_{(Z)}(V). 
$$
\end{proposition}

%---------------------------------------
\begin{lemma} \label{lemma-3}  %lemma4.4   
%---------------------------------------
Let $A=\oplus_{i=0}^c A_i$ be a graded Artinian $K$-algebra 
with the WLP (resp. SLP) 
and let $y, z \in A_1$ be two linear forms of $A$. 
If $y$ is a weak Lefschetz element (resp. strong Lefschetz element) for  $A$, 
then so is $y+\lambda z$ for a general  element $\lambda\in K$.  
\end{lemma}

\begin{proof}
Let $M$ and $N$ be the square matrices of 
$\times y \in \End(A)$ and $\times z \in \End(A)$ 
for a same basis of $A$, respectively. 
Suppose that $y$ is a strong Lefschetz element for  $A$. 
By assumption and Lemma~\ref{lemma2-2}, 
we have that $\rank(M^k)=\dim A-{\rm SP}_k(A)$ 
for all $k=1,2,\ldots,c$. 
We would like to show that 
$\rank (M+\lambda N)^k=\dim A-{\rm SP}_k(A)$ 
for a general element $\lambda\in K$ and all $k=1,2,\ldots,c$. 
Let $Q_0$ be an $n\times n$ square submatrix of $M^k$ 
such that $\rank(Q_0)=\rank M^k$ and $Q_0$ is regular. 
Put 
$$
P=(M+\lambda N)^k=M^k+M^{k-1}N\lambda+\cdots+MN^{k-1}\lambda^{k-1}+N^k\lambda^k,
$$
and let $P^\prime$ and $Q_{i}$ be the submatrices of $P$ and $M^{k-i}N^i$ $(1\leq i\leq k)$, 
respectively, 
obtained by deleting the same rows and colums as $Q_0$. 
Obviously 
$$
P^\prime=Q_0+Q_1\lambda+\cdots+Q_k\lambda^k. 
$$
Hence, we have that 
$$
\det(P^\prime) = 
\alpha_0 + \alpha_1\lambda + \alpha_2\lambda^2 + \cdots + \alpha_{nk}\lambda^{nk} 
$$
for some $\alpha_i\in K$ $(0\leq i\leq nk)$, 
where $\alpha_0=\det(Q_0)$. 
Thus, noting that $\det(Q_0)\neq 0$, 
our assertion follows from Lemmas~\ref{lemma2-2} and~\ref{lemma2.5}. 

Next suppose that $y$ is a weak Lefschetz element for  $A$. 
By assumption and Lemma~\ref{lemma2-1}, 
we have that $\rank(M)={\rm CoSperner}(A)$. 
We would like to show that 
$\rank (M+\lambda N)={\rm CoSperner}(A)$ 
for a general element $\lambda\in K$. 
This also follows by the same idea as above. 
\end{proof}

\begin{notation and remark} \label{notation and remark 3-3}  %remark4.5  
Let $A=\oplus_{i=0}^c A_i$ be a graded Artinian $K$-algebra, 
let $z\in A_1$ be a linear form of $A$ 
and let $P(\times z)=n_1\oplus n_2\oplus\cdots\oplus n_r$ 
be the Jordan decomposition of $\times z \in \End(A)$. 
Now we can take homogeneous elements 
$a_i\in A\setminus (z)$ $(i=1,2,\ldots,r)$ such that 
the set 
\begin{equation} \label{EQ:3-7-1}
\cup_{i=1}^r\{a_i,a_iz,a_iz^2,\ldots,a_iz^{n_i-1}\}
\end{equation}
is a basis for $A$ as a vector space, that is, 
the matrix of $\times z \in \End(A)$ for the basis above 
coincides with the Jordan canonical from of $\times z$. 
Then it is easy to check that the set 
$$
\cup_{i=1}^r\{a_i^\ast,a_i^\ast z^\ast,a_i^\ast(z^\ast)^2,\ldots,
a_i^\ast(z^\ast)^{n_i-1}\}
$$
is a basis for $\Gr_{(z)}(A)$. 
Hence the Jordan canonical form of $\times z^{\ast} \in \End(       \Gr_{(z)}(A)        )$ is 
the same as that of $\times z \in \End(A)$. 
\end{notation and remark}

%%%%%%%%%%%%%%%%%%%%%%%%%%%%%%%%%%%%%%%%%%%%%%%%%%%%%%%%%%%%%%%%%%%%%%%%%%%%%%
%%%%%%%%%%%%%%%%%%%%%%%%%%%%%%%%%%%%%%%%%%%%%%%%%%%%%%%%%%%%%%%%%%%%%%%%%%%%%%
%---------------------------------------
\begin{theorem} \label{main-th11}  %theorem1.2   %thm4.6 
Let $K$ be a field of characteristic zero 
and let $A$ be a graded Artinian $K$-algebra. 
Then 
$A$ has the WLP (resp. SLP) if and only if 
$\Gr_{(z)}(A)$ has the WLP (resp. SLP) for some linear form $z$ of $A$. 
\end{theorem}
%---------------------------------------
%%%%%%%%%%%%%%%%%%%%%%%%%%%%%%%%%%%%%%%%%%%%%%%%%%%%%%%%%%%%%%%%%%%%%%%%%%%%%%
%%%%%%%%%%%%%%%%%%%%%%%%%%%%%%%%%%%%%%%%%%%%%%%%%%%%%%%%%%%%%%%%%%%%%%%%%%%%%%
\begin{proof}  %mainthm_1
($\Rightarrow$)\ 
Let $z$ be a  Lefschetz element of $A$. 
Recall that  $A$ and $\Gr_{(z)}(A)$ have the same Hilbert function. 
Then Remark~\ref{notation and remark 3-3} and Lemma~\ref{lemma3} shows that 
$z^{\ast}$ is a Lefschetz element of $\Gr_{(z)}(A)$.

($\Leftarrow$)\ 
{\em Step 1: }  
First we prove that if $\Gr_{(z)}(A)$ has the WLP then so does $A$. 
Let $g$ be a weak Lefschetz element of $\Gr_{(z)}(A)$. 
Since $[\Gr_{(z)}(A)]_1=A_1/zA_0\oplus zA_0$, 
$g$ can be written as $g=y^\ast+\lambda_0 z^\ast$ 
for some $y\in A_1$ and $\lambda_0 \in K$. 
Let $S=K[Y,Z]$ be the polynomial ring in two variables. 
Define the  algebra homomorphism $S \ra \End(A)$  
by $Y \mapsto \times y$ and $Z \mapsto \times z$. 
Then we may regard $A$ as a finite graded $S$-module. 
>From Proposition~\ref{prop}, 
it follows that 
\begin{equation}\label{EQ}
\dim A/(y+\lambda z)A \leq 
\dim \Gr_{(z)}(A)/(y^{\ast}+\lambda z^{\ast}) \Gr_{(z)}(A)
\end{equation}
for a general element $\lambda \in K$. 
Furthermore, 
it follows from Lemmas~\ref{lemma2-1}  and \ref{lemma-3} that 
$$
\dim \Gr_{(z)}(A)/(y^{\ast}+\lambda z^{\ast}) \Gr_{(z)}(A) = 
{\rm Sperner}(\Gr_{(z)}(A)) 
$$ 
for a general element $\lambda\in K$. 

On the other hand, since 
${\rm Sperner}(\Gr_{(z)}(A)) = {\rm Sperner}(A)$ and since 
${\rm Sperner}(A) \leq \dim A/(y + \lambda z)A$ by Lemma~\ref{lemma2.5}, 
we have $${\rm Sperner}(A) = \dim A/(y + \lambda z)A,$$
proving that $A$ has the WLP.  
\medskip

\noi
{\em Step 2:} 
It still remains to show that $A$  has the SLP 
assuming that $\Gr_{(z)}(A)$ has the SLP. 
Let $t$ be a new variable and let $\widetilde{A}=A[t]/(t^{\alpha})$, 
where $\alpha$ is any positive integer. 
Then, since we have 
$$\Gr_{(z)}(\widetilde{A})\cong \Gr_{(z)}(A) \otimes _K K[t]/(t^{\alpha})$$
and since the SLP is preserved by tensor product,   
it follows that $\Gr_{(z)}(\widetilde{A})$ has the SLP. 
By Step 1, 
this implies that 
$\widetilde{A}$ has the WLP for all $\alpha > 0$.
Hence the SLP of $A$ follows by Proposition~\ref{lemma_SLP_WLP}.  
\end{proof}

%%%%%%%%%%%%%%%%%%%%%%%%%%%%%%%%%%%%%%%%%%%%%%%%%%%%%%%%%%%%%%%%%%%
%
%
%%%%    End of Section 4  --- start of Section 5 %%%%%%  
%
%
%%%%%%%%%%%%%%%%%%%%%%%%%%%%%%%%%%%%%%%%%%%%%%%%%%%%%%%%%%%%%%%%%%%

\section{Central simple modules and the strong Lefshetz property}  %%%section5
In this section we prove Theorems~\ref{main-th2} and \ref{new_thm}.  

\begin{definition and remark} \label{Jordan sequence}   %notation5.1
Let $A=\oplus_{i=0}^c A_i$ be a graded Artinian $K$-algebra, 
let $z\in A_1$ be a linear form of $A$ 
and let $P(\times z)=n_1\oplus n_2\oplus\cdots\oplus n_r$ 
be the Jordan decomposition of $\times z \in \End(A)$. 
Furthermore, 
let $(f_1,f_2,\ldots,f_s)$ 
be the finest subsequence of $(n_1,n_2,\ldots,n_r)$ 
such that $f_1>f_2>\cdots>f_s$. 
Then we rewrite the same Jordan decomposition $P(\times z)$ as
$$
%\begin{equation}   \label{grand_block_decomposition}   %Equation10
P(\times z)=
n_1\oplus n_2\oplus\cdots\oplus n_r=
\underbrace{f_1\oplus\cdots\oplus f_1}_{m_1}\oplus 
\underbrace{f_2\oplus\cdots\oplus f_2}_{m_2}\oplus\cdots\oplus
\underbrace{f_s\oplus\cdots\oplus f_s}_{m_s}. 
%\end{equation}
$$
We call the graded $A$-module 
$$
U_i=\displaystyle\frac{(0:z^{f_i})+(z)}{(0:z^{f_{i+1}})+(z)} 
$$ 
the {\em $i$-th central simple module} of $(A,z)$, 
with $1\leq i\leq s$ and $f_{s+1}=0$.  
Note that  these are defined for a pair of the algebra $A$ 
and a linear form $z\in A_1$. 
By the definition, 
it is easy to see that 
the modules $U_1,U_2,\ldots,U_s$ are the non-zero terms of 
the successive quotients of the descending chain of ideals 
$$
A=(0:z^{f_1})+(z) \supset (0:z^{f_1-1})+(z) \supset \cdots 
\supset (0:z)+(z) \supset (z). 
$$

For $1 \leq i \leq s$, define $\widetilde{U_i}$  by   
$\widetilde{U_i}=U_i\otimes_K K[t]/(t^{f_i}).$  
\end{definition and remark}

%%%%%%%%%%%%%%%%%%%%%%%%%%%%%%%%%%%%%%%%%%%%%%%%%%%%%%%%%%%%%%%%%%%%%%%%%%%%%A
%%%%%%%%%%%%%%%%%%%%%%%%%%%%%%%%%%%%%%%%%%%%%%%%%%%%%%%%%%%%%%%%%%%%%%%%%%%%%A
%--------------------------------------
\begin{theorem} \label{main-th2}   %theorem1.2  %mainthm1.2 thm5.2
Let $K$ be a field of characteristic zero 
and let $A$ be a graded Artinian $K$-algebra. 
Then the following conditions are equivalent.  
\begin{itemize}
\item[(i)] 
$A$ has the SLP. 
\item[(ii)]
There exists a linear form $z$ of $A$ such that 
all the central simple modules $U_i$ of $(A,z)$ have the SLP  and 
the reflecting degree of the Hilbert function of $\widetilde{U_i}$ 
coincides with that of $A$ for $i= 1, 2, \ldots, s$.   
%where $P(\times z)=f_1\oplus\cdots\oplus f_1\oplus 
%f_2\oplus\cdots\oplus f_2\oplus\cdots\oplus f_s\oplus\cdots\oplus f_s$ 
%is the Jordan decomposition of $\times z:A \ra A$ 
%and put $\widetilde{U_i}=U_i\otimes_K K[t]/(t^{f_i})$ 
%XSfor all $i=1,2,\ldots,s$. 
\end{itemize}
\end{theorem}
%--------------------------------------
%%%%%%%%%%%%%%%%%%%%%%%%%%%%%%%%%%%%%%%%%%%%%%%%%%%%%%%%%%%%%%%%%%%%%%%%%%%%%A

\begin{proof}  %[Outline of proof of Theorem~\ref{main-th2}] 
(i)$\Rightarrow$(ii): 
Assume that $A$ has the SLP and $z$ is a strong Lefschetz element. 
Put 
$\overline{A}_i =A_i/zA_{i-1}$.  
Then we may write 
$$A/(z)= \oplus _{i=0}^ {c'}\overline{A_i}$$
where $c'$ is the largest integer such that 
$(A/(z))_{c'} \not = 0$. 
Then, noting that $A$ has the SLP, 
one sees easily that $(A,z)$ has  
$c'+1$ central simple modules $U_1, \ldots, U_{c'+1}$ 
and that  $U_i \cong \ol{A_{i-1}}$.  
This shows that $U_i$ has only one non-trivial 
graded piece concentrated at the degree $i-1$. 
Hence $U_i$ has the SLP for trivial reasons. 
Also it is easy to show that 
the reflecting degree of the  Hilbert function of $\widetilde{U_i}$ 
coincides with that of $A$ 
for all $i=1,2,\ldots,c'+1$. 
\medskip

\noi
(ii)$\Rightarrow$(i): By  Theorem~\ref{main-th11} it suffices to prove that $\Gr_{(z)}(A)$ has the SLP.  
%We use Notation~\ref{Jordan sequence}. 
We divide the proof into three steps. 
\medskip

\noi
{\em Step1}: 
Choose $g \in A_1$ general enough so that 
$g$ is a strong Lefschetz element of $U_i$ for any  $i=1,2,\ldots,s$.   
Let $m_i=\dim U_i$ $(1\leq i\leq s)$. 
Now we take a basis of $U_i$ as a vector space for all $1\leq i\leq s$, 
$$
\{ \ol{e_{i1}}, \ol{e_{i2}}, \ldots, \ol{e_{im_{i}}} \}, 
$$
where $e_{ij}\in (0:z^{f_i})+(z)$, $j=1, \ldots, m_i$.  
Then the set 
$$
\{ \ol{e_{ij}}\otimes \overline{t}^k \mid 1\leq j\leq m_i, \ 0\leq k\leq f_i-1\}
$$
is a basis of $\ti{U_i}=U_i\otimes_K K[t]/(t^{f_i})$ 
for all $i$ such that  $1\leq i\leq s$. 
Put $\ti{U}=\oplus_{i=1}^s\ti{U_i}$. 
The set 
\begin{equation} \label{basis1}
\cup_{i=1}^s\{ \ol{e_{ij}}\otimes \overline{t}^k \mid 1\leq j\leq m_i, \ 0\leq k\leq f_i-1\}
\end{equation}
is a basis of $\ti{U}$. 
We may consider $\ti{U}$ as a graded module over  $A\otimes K[t]$. 
Here we calculate a matrix of the multiplication 
$\times (g\otimes 1 + 1\otimes t) :  \ti{U} \ra \ti{U}$ 
as  an endomorphism.  
Let $P_i$ be the square matrix of 
$\times (g\otimes 1 + 1\otimes t) : \ti{U_i}\ra\ti{U_i}$ 
using  the basis above. 
Since $(g\otimes 1 + 1\otimes t)\ti{U_i}\subset\ti{U_i}$, 
a matrix for $\times (g\otimes 1 + 1\otimes t) \in \End(\ti{U})$ 
is of the following form,  
$$
P=
\left[
\begin{array}{ccccccc}
P_1 &     &        & \bigzerou  \\
    & P_2 &        &            \\
    &     & \ddots &            \\
\bigzerol   &      &        & P_s   \\     
\end{array}
\right]. 
$$
Hence it follows that 
$$
P^h=
\left[
\begin{array}{ccccccc}
P_1^h &     &        & \bigzerou  \\
    & P_2^h &        &            \\
    &     & \ddots &            \\
\bigzerol   &      &        & P_s^h   \\     
\end{array}
\right]
$$
for all $h$. 
\medskip

\noi
{\em Step 2}: 
Let $g^\ast$ be the initial form  of $g$  in $\Gr_{(z)}(A)$. 
We calculate a matrix of the multiplication 
$\times g^\ast: \Gr_{(z)}(A)\ra \Gr_{(z)}(A)$. 
First we note that the set 
\begin{equation} \label{basis2}
\cup_{i=1}^s \{ (e_{ij}^\ast)(z^\ast)^k \mid 1\leq j\leq m_i, 
\ 0\leq k\leq f_i-1 \}
\end{equation}
is a basis of $\Gr_{(z)}(A)$. 
Let $V_i$ be the  subspace of $A/(z)$ in $\Gr_{(z)}(A)$ 
which is generated by $\{e_{i1}^\ast, e_{i2}^\ast, \ldots, e_{im_i}^\ast \}$ 
for all $i=1,2,\ldots,s$. 
Furthermore, 
let $V_i^\ast$ be the subspace of $\Gr_{(z)}(A)$ which is generated by 
$$
\{ (e_{ij}^\ast)(z^\ast)^k \mid 1\leq j\leq m_i, \ 0\leq k\leq f_i-1 \}
$$ 
for all $i=1,2,\ldots,s$. 
Then, since $g^\ast V_i\subset\oplus_{j=i}^s V_j$ for all $i=1,2,\ldots,s$, 
we have that $g^\ast V_i^\ast\subset\oplus_{j=i}^s V_j^\ast$.  
Hence, a matrix of $\times (g^\ast + z^\ast) : \Gr_{(z)}(A)\ra \Gr_{(z)}(A)$ is 
of the following form: 
$$
Q=
\left[
\begin{array}{ccccccc}
P_1 &     &        & \bigast  \\
    & P_2 &        &            \\
    &     & \ddots &            \\
\bigzerol   &      &        & P_s   \\     
\end{array}
\right]. 
$$
Thus, we obtain that 
$$
Q^h=
\left[
\begin{array}{ccccccc}
P_1^h &     &        & \bigast  \\
    & P_2^h &        &            \\
    &     & \ddots &            \\
\bigzerol   &      &        & P_s^h   \\     
\end{array}
\right]
$$
for all $h$. 
\medskip

\noi
{\em Step 3}: 
Note that  
$\ti{U}$ and $\Gr_{(z)}(A)$ have the same Hilbert function, since 
$\deg \ \ol{e_{ij}}=\deg \ e_{ij}^\ast$ for all $i$ and $j$ 
and  since the sets (\ref{basis1}) and (\ref{basis2}) 
are bases for $\ti{U}$ and $\Gr_{(z)}(A)$ respectively. 

Next we claim that  
$\ti{U} = \oplus \ti{U}_i$ has the SLP as an $A\otimes_KK[t]$-module. 
Indeed, by Proposition~\ref{proposition3.4}, every 
$\ti{U_i}$ has the SLP as an $A\otimes_KK[t]$-module 
with  $g\otimes 1+1\otimes t$  a strong Lefschetz element. 
By assumption 
the  reflecting degree of $h_{\ti{U_i}}(q)$ 
coincides  with that of $h_{\ti{U}}(q)$ for any $i$.  
This proves that the element 
$g\otimes 1+1\otimes t$ is a strong Lefschetz element for  $\ti{U}=\oplus \ti{U}_i$. 

Now if $\rank \ P^h \leq \rank \ Q^h$ for all $h$,  it implies that $\Gr_{(z)}(A)$ has the SLP by 
Lemma~\ref{lemma2.5}~(4).  So we prove that   $\rank \ P^h \leq \rank \ Q^h$ for all $h$.  

Let $P_{i h}^\prime$ be a square submatrix of $P_i^h$ 
such that $P_{i h}^\prime$ is of full rank  and 
$\rank \ P_i^h = \rank \ P_{i h}^\prime$ 
for all $i$ and $h$. 
Then we have  
$$
\begin{array}{rcl}
\rank \ P^h & = & \sum_{i=1}^s \rank \ P_{i}^h \\[1ex]
& = & \sum_{i=1}^s \rank \ P_{i h}^\prime \\[1ex]
& = & \rank \ P_h^\prime,  
\end{array}
$$ 
where 
$$
P_h^\prime=
\left[
\begin{array}{ccccccc}
P_{1 h}^\prime &     &        & \bigzerou  \\
    & P_{2 h}^\prime &        &            \\
    &     & \ddots &            \\
\bigzerol   &      &        & P_{s h}^\prime   \\     
\end{array}
\right]. 
$$
Let $Q_h^\prime$ be the square submatrix of $Q^h$ 
consisting of the same rows and columns as $P_h^\prime$, 
so that  
$Q_h^\prime$ is of  the following form 
$$
Q_h^\prime=
\left[
\begin{array}{ccccccc}
P_{1 h}^\prime &     &        & \bigast  \\
    & P_{2 h}^\prime &        &            \\
    &     & \ddots &            \\
\bigzerol   &      &        & P_{s h}^\prime   \\     
\end{array}
\right]. 
$$
Then, since $\det Q_h^\prime =\prod_{i=1}^s \det P_{i h}^\prime \neq 0$, 
it follows that $\rank \ Q_h^\prime = \rank \ P_h^\prime = \rank \ P^h$. 
This means that $\rank \ P^h \leq \rank \ Q^h$ for all $h$. 
\end{proof}

\begin{proposition}  \label{csm_of_Gorenstein_algebra}  %prop5.3
Let $K$ be any field and suppose that $A$ is an graded Artinian Gorenstein $K$ algebra. 
Suppose that $z$ is a linear form of $A$ and 
let  
$$
%\begin{equation}   \label{grand_block_decomposition}   %Equation10
P(\times z)=
%n_1\oplus n_2\oplus\cdots\oplus n_r=
\underbrace{f_1\oplus\cdots\oplus f_1}_{m_1}\oplus 
\underbrace{f_2\oplus\cdots\oplus f_2}_{m_2}\oplus\cdots\oplus
\underbrace{f_s\oplus\cdots\oplus f_s}_{m_s}. 
%\end{equation}
$$
be as in  
{\rm Notation~\ref{Jordan sequence}}. 
Furthermore let $U_i$ be the $i$th central simple module of $(A,z)$.  
Then, for any $i$ with $1 \leq i \leq s$, the graded module  $U_i$ has a symmetric Hilbert function and  
the module $\ti{U}_i=U_{i} \otimes K[t]/(t^{f_i})$ has a symmetric Hilbert function with the 
same reflecting degree as that of $A$. 
\end{proposition} 

\begin{proof}
Let $c$ be the socle degree of $A$.  
We induct on $c$.    
 First note that  $A/0:z$ is Gorenstein and  
the socle degree of $A/0:z$  is   $c-1$, as long as  $A/0:z \neq 0$.  
To see this consider the exact sequence  
\[
0 \ra (A/0:z)(-1) \ra A \ra A/zA \ra 0, 
\]
where the first map sends $1$ to $z$.  
 Since the zero ideal of  $A$ is irreducible, 
the sequence shows that the zero ideal of $A/0:z$ is irreducible also. 
This implies  that   $A/0:z$   is a Gorenstein algebra (cf. \cite{jW73}).   
Recall that $A_c$, the socle,  is the unique minimal ideal of $A$.  
Hence  the  ideal $zA \cong (A/0:z)(-1)$  contains $A_c$, from 
which it follows that  $(A/zA)_{c}=0$.  This shows that $zA_{c-1} \neq 0$.  On the other hand, obviously, it holds 
that $zA_{c}=0$.  Thus  the socle degree of $A/0:z$ is precisely $c-1$.  

If $A/0:z = 0$, then $s=1$, $f_1=1$ and  $U_1 = \ti{U_1} =A$, and the assertions  
are  clear from the fact that $A$ has a symmetric Hilbert function.  
Now  we assume that $A/0:z \neq 0$.  
Note that the Jordan-block-size for $\times \ol{z} \in \End(A/0:z)$ is given by 
$$
P(\times \ol{z})=
%n_1\oplus n_2\oplus\cdots\oplus n_r=
\underbrace{f_1-1\oplus\cdots\oplus f_1-1}_{m_1}\oplus 
\underbrace{f_2-1\oplus\cdots\oplus f_2-1}_{m_2}\oplus\cdots\oplus
\underbrace{f_s-1\oplus\cdots\oplus f_s-1}_{m_s}. 
$$
First assume that $f_s >1$.   
Then  $(A/0:z, \ol{z})$ and $(A,z)$ have the same central simple modules.  Thus 
$U_i$ has the symmetric Hilbert function for all $i$ such that $1 \leq i \leq s$.  
Moreover the reflecting  degrees of $A/0:z$ and $A$ differ by $1/2$. Hence the second assertion  is also clear. 
Now assume that $f_s=1$.  It is possible to regard  $U_1, \ldots, U_{s-1}$  as the central simple modules for 
$(A/0:z, \ol{z})$.  Thus, by induction hypothesis,  $U_i \otimes K[t]/(t^{f_i-1})$,  for all $i=1,2, \ldots, s-1$,  
have the same reflecting degree as that of 
$A/0:z$.  
Hence the same is true  for $A$ and  $\ti{U}_i$  with  $i=1,2, \ldots, s-1$.  
It remains to prove that, with the assumption $f_s=1$, the modules $U_s$ and   $\ti{U_s}$ 
have symmetric Hilbert functions, and  $\ti{U_s}$   with the same reflecting degree as that of $A$.  
   Since $f_s=1$, we have $\ti{U}_s=U_s$.   
Notice that 
$$
h_A(q)=\sum _{j=1} ^ {s} h_{\ti{U}_j}(q).
$$
We may apply the induction hypothesis to the first $s-1$ summands.  Hence the last summand 
$h_{\ti{U}_s}(q)$ also has the reflecting degree equal to that of $A$. 
\end{proof}

\begin{theorem} \label{new_thm}   %theorem5.4  %thm5.4 
Let $K$ be a field of characteristic zero 
and let $A$ be a graded Artinian Gorenstein $K$-algebra. 
Then the following conditions are equivalent.  
\begin{itemize}
\item[(i)] 
$A$ has the SLP. 
\item[(ii)]
There exists a linear form $z$ of $A$ such that 
all the central simple modules $U_i$ of $(A,z)$ have the SLP.
%XSfor all $i=1,2,\ldots,s$. 
\end{itemize}
\end{theorem}

\begin{proof}
This follows from Theorem~\ref{main-th2} and Proposition~\ref{csm_of_Gorenstein_algebra}. 
\end{proof}

%%%%%%%%%%%%%%%%%%%%%%%%%%%%%%%%%%%%%%%%%%%%%%%%%%%%%%%%%%%%%%%%%%%
%
%
%%%%    End of Proof of Main Theorem  --- start of Section 6 %%%%%%
%
%
%%%%%%%%%%%%%%%%%%%%%%%%%%%%%%%%%%%%%%%%%%%%%%%%%%%%%%%%%%%%%%%%%%%

\section{finite free extensions of a graded Artinian $K$-algebra}  %section6 

Theorem~\ref{main-th3} below is an extension of Theorem 28 in \cite{tHjW04}. 
We can now give another  proof simplifying very much the proof given in  \cite{tHjW03} and  \cite{tHjW04}.

%mainthm1.3
%%%%%%%%%%%%%%%%%%%%%%%%%%%%%%%%%%%%%%%%%%%%%%%%%%%%%%%%%%%%%%%%%%%%%%%%
%%%%%%%%%%%%%%%%%%%%%%%%%%%%%%%%%%%%%%%%%%%%%%%%%%%%%%%%%%%%%%%%%%%%%%%%
%--------------------------------------
\begin{theorem} \label{main-th3}   %theorem1.3   %thm6.1
Let $K$ be a field of characteristic zero, 
let $B$ be a graded Artinian $K$-algebra 
and let $A$ be a finite flat algebra over $B$ 
such that the algebra map $\varphi:B \ra A$ preserves grading. 
Assume that both $B$ and $A/mA$ have the SLP, 
where $m$ is the maximal ideal of $B$. 
Then $A$ has the SLP. 
\end{theorem}
%-------------------------------------
%%%%%%%%%%%%%%%%%%%%%%%%%%%%%%%%%%%%%%%%%%%%%%%%%%%%%%%%%%%%%%%%%%%%%%%%
%%%%%%%%%%%%%%%%%%%%%%%%%%%%%%%%%%%%%%%%%%%%%%%%%%%%%%%%%%%%%%%%%%%%%%%%

First we prove a lemma.   

\begin{lemma} \label{lemma6-1}  %lemma6.2
We use the same notation as {\rm Theorem~\ref{main-th3}}. 
Let $z'$ be any linear form of $B$ 
and put $z=\varphi(z')$. 
Let $U_i'$ and $U_i$ be 
the $i$-th central simple modules of 
$(B,z')$ and $(A,z)$, respectively. 
Then $U_i'\otimes_BA \cong U_i$. 
\end{lemma}

\begin{proof}
By assumption, 
we can write 
$A\cong Be_1\oplus Be_2\oplus\cdots\oplus Be_k$ 
for some homogeneous elements $e_i\in A$. 
Let 
$$
%\begin{equation}   \label{grand_block_decomposition}   %Equation10
P(\times z')=
\underbrace{f_1\oplus\cdots\oplus f_1}_{m_1}\oplus 
\underbrace{f_2\oplus\cdots\oplus f_2}_{m_2}\oplus\cdots\oplus
\underbrace{f_s\oplus\cdots\oplus f_s}_{m_s}
%\end{equation}
$$
be the Jordan decomposition of $\times z' \in \End(B)$. 
Then it follows immediately that 
$$
%\begin{equation}   \label{grand_block_decomposition}   %Equation10
P(\times z)=
\underbrace{f_1\oplus\cdots\oplus f_1}_{m_1\times k}\oplus 
\underbrace{f_2\oplus\cdots\oplus f_2}_{m_2\times k}\oplus\cdots\oplus
\underbrace{f_s\oplus\cdots\oplus f_s}_{m_s\times k}. 
%\end{equation}
$$
Hence, noting that 
$$
(0:_Az^m)+zA=\{(0:_B(z')^m)+z'B\}e_1
\oplus\cdots\oplus\{(0:_B(z')^m)+z'B\}e_k 
$$
for all integers $m>0$, 
we have that 
$$
\begin{array}{rcl}
U_i'\otimes_BA & \cong & 
[\{(0:_B(z')^{f_i})+z'B\}/\{(0:_B(z')^{f_{i+1}})+z'B\}]
\otimes_B(\oplus_{i=1}^k Be_i) \\[1ex]
    & \cong &  
\oplus_{i=1}^k
[\{(0:_B(z')^{f_i})+z'B\}/\{(0:_B(z')^{f_{i+1}})+z'B\}]e_i \\[1ex]
                & \cong & 
\{(0:_Az^{f_i})+zA\}/\{(0:_Az^{f_{i+1}})+zA\} \\[1ex]
  & \cong  & U_i              
\end{array}
$$
for all $i=1,2,\ldots,s$. 
\end{proof}

\begin{proof}
[Proof of Theorem~\ref{main-th3}] 
Let $z'$ be a strong Lefschetz element of $B$ 
and put $z=\varphi(z')$. 
>From the proof (i)$\Rightarrow$(ii) of Theorem~\ref{main-th2},  
every central simple module $U_i'$ of $(B,z')$ has 
only one non-trivial graded piece. 
Hence $U_i'$ has the SLP. 
Also, since $U_i'$ is annihilated by $m$, 
we have by Lemma~\ref{lemma6-1} that $U_i'\otimes_KA/mA \cong U_i$, 
where $U_i$ is the $i$-th central simple module of $(A,z)$. 
Thus, using our assumption that $A/mA$ has the SLP, 
it follows by Proposition~\ref{proposition3.4} 
that every central simple module of $(A,z)$ has the SLP. 

To proceed with the proof, 
we use the same notation as the proof of Lemma~\ref{lemma6-1}. 
Fix $i$ such that  $1 \leq i \leq s$ and put 
$$
\ti{U_i}=U_i\otimes_KK[t]/(t^{f_i}) 
\ \ {\rm and}\ \ 
\ti{U_i'}=U_i'\otimes_KK[t]/(t^{f_i}).  
$$
Let $a,b$ be the initial and end degrees of the Hilbert function of $\ti{U}'_i$ so that  
$$h_{\ti{U_i'}}(q)=h_aq^a   +   \mbox{ \ (\rm mid terms) \ } +  h_bq^b.$$ 
Similarly let  
$$
h_B(q)= 1+ \mbox{\rm (mid terms) } + q^c, 
\ \ 
h_{A/mA}(q)= 1+ \mbox{\rm (mid terms) } + q^d
$$
be the Hilbert functions of $B$ and $A/mA$. 
Since  
$h_{\ti{U_i}}(q)=h_{\ti{U_i'}}(q)h_{A/mA}(q)$ 
and $h_A(q)=h_B(q)h_{A/mA}(q)$, 
the Hilbert functions of $\ti{U_i}$ and $A$ are of  the following form:  

$$
h_{\ti{U_i}}(q)=h_a q^a+   \mbox{ \ (\rm mid terms) \ }   + h_bq^{b+d},   
$$
$$ 
h_A(q)=1 +  \mbox{ \ (\rm mid terms) \ }   + q^{c+d}. 
$$  
%$$
%h_{\ti{U_i}}(q)=h_a q^a+\cdots+(h_b+\beta)q^{b+d}  
%\ \ {\rm and}\ \ 
%h_A(q)=(\alpha+\beta)q^{c+d}+{\rm (lower \ terms)}. 
%$$  
Since the reflecting degree of $h_{\ti{U_i'}}(q)$ 
coincides with that of $h_B(q)$, 
we have  $(a+b)/2=c/2$. 
Hence we obtain $(a+b+d)/2=(c+d)/2$. 
This means that 
the reflecting degree of $h_{\ti{U_i}}(q)$ 
coincides with that of $h_A(q)$ independent of $i$.   
Thus by Theorem~\ref{main-th2}  proof is finished. 
\end{proof}

\begin{example}   %ex6.3 
In general the converse of Theorem~\ref{main-th3} is not true. 
We give such an example.  
Let $A=K[x_1,x_2]/(e_1^2,e_2^2)$ and $B=K[e_1,e_2]/(e_1^2,e_2^2)$, 
where $e_1=x_1+x_2$ and $e_2=x_1x_2$. 
Note that $A$ has the standard grading, but $B$ does not.  
By Lemma~\ref{lemma:7-4} below, 
$A$ is a finite free module over $B$. 
Furthermore, by Propositon 4.4 of \cite{tHjMuNjW01}, 
both $A$ and $A/mA$ have the SLP. 
However, $B$ does not have the SLP.   
In fact $e_1$ is the only candidate for a strong Lefschetz element for 
$B$ while we have  $h_B(q)=1+q+q^2+q^3$ and $\overline{e_1}^2=0$ in $B$.  
\end{example}

\begin{remark}   %rem6.4  
The converse of Theorem~\ref{main-th3} is true 
under some assumptions. 
We use the same notation as Theorem~\ref{main-th3}. 
Assume that there is a linear form $z'$ of $B$ such that 
\begin{itemize}
\item[(i)]
all central simple modules of ($A,\varphi(z')$) 
have the SLP, 
\item[(ii)]
every central simple module of ($B,z'$) has 
a symmetric unimodal Hilbert function and 
\item[(iii)]
$A/mA$ has a symmetric unimodal Hilbert function. 
\end{itemize}
Then, 
using Theorems~\ref{main-th2} and~\ref{main_thm_one} 
and the same idea as the proof of Theorem~\ref{main-th3}, 
we can easily prove that $B$ and $A/mA$ have the SLP. 
\end{remark}

%%%%%%%%%%%%%%%%%%%%%%%%%%%%%%%%%%%%%%%%%%%%%%%%%%%%%%%%%%%%%%%%%%%
%
%
%%%%    End of Section 6  --- start of Section 7 %%%%%%
%
%
%%%%%%%%%%%%%%%%%%%%%%%%%%%%%%%%%%%%%%%%%%%%%%%%%%%%%%%%%%%%%%%%%%%

\section{Complete intersections defined by power sums 
of consecutive degrees}

Using Theorem~\ref{main-th3} we prove the following: 

%prop1.4
%%%%%%%%%%%%%%%%%%%%%%%%%%%%%%%%%%%%%%%%%%%%%%%%%%%%%%%%%%%%%%%%%%%%%%%
%%%%%%%%%%%%%%%%%%%%%%%%%%%%%%%%%%%%%%%%%%%%%%%%%%%%%%%%%%%%%%%%%%%%%%%
%-------------------------------------
\begin{proposition} \label{main-th4}   %proposition1.4  %prop7.1
Let $R=K[x_1,x_2,\ldots,x_n]$ be the polynomial ring 
over a field $K$ of characteristic zero 
with $\deg(x_i)=1$ for all $i$. 
Let $a$ be a positive integer. 
Let
$$
I=(p_a(x_1,\ldots,x_n), p_{a+1}(x_1,\ldots,x_n), 
\ldots, p_{a+n-1}(x_1,\ldots,x_n)), 
$$
where $p_d=x_1^d+x_2^d+\cdots+x_n^d$ is the power sum symmetric function  of degree $d$. 
Then  $A=R/I$ is a complete intersection and has the SLP.  
\end{proposition}
%-------------------------------------
%%%%%%%%%%%%%%%%%%%%%%%%%%%%%%%%%%%%%%%%%%%%%%%%%%%%%%%%%%%%%%%%%%%%%%%
%%%%%%%%%%%%%%%%%%%%%%%%%%%%%%%%%%%%%%%%%%%%%%%%%%%%%%%%%%%%%%%%%%%%%%%

The proof is given at the end of this section after a series of lemmas.

\begin{lemma} \label{lemma:7-1}   %lemma7.2
With the same noation as  {\rm  Proposition~\ref{main-th4}},  
 the ideal $I$ is a complete intersection. 
\end{lemma}

\begin{proof}
The well known identity

\[
p_m=-\sum _{j=1} ^n (-1)^j e_j p_{m-j}
\] 
where  $e_i$ is the elementary symmetric function of degree $i$,   
implies that 
$p_m\in I$ for all $m>a+n-1$. 
Hence we have 
$(p_m, p_{2m}, \ldots, p_{nm})\subset I$. 
It is easy to see  that 
 $(p_m, p_{2m}, \ldots, p_{nm})$ is  a complete intersection and hence so is $I$. 
\end{proof}

\begin{lemma} \label{lemma:7-2}   %lemma7.3
Let $C=\oplus_{i=0}^cC_i$ be a graded Artinian $K$-algebra 
with the SLP, where $C_c\neq (0)$. 
Let $B=\oplus_{i=0}^cB_i$ be a graded subalgebra of $C$, where $B_c\neq (0)$. 
If  $B$ contains a strong Lefschetz element for $C$ and if the 
Hilbert function of $B$ is symmetric, then $B$ has the SLP. 
\end{lemma}

\begin{proof}
Suppose $z \in B$ is a strong Lefschetz element for $C$. 
Then 
the multiplication $\times z^{c-2i}: B_i\ra B_{c-i}$ 
is injective for all  $i \leq [c/2]$. 
Hence, noting that $\dim B_i=\dim B_{c-i}$, 
we have that 
$\times z^{c-2i}$ is bijective for all $i$. 
\end{proof}

\begin{notation and remark} \label{NR:7-3}  %remark7.4
Let $K$ be a field of characteristic zero 
and $R=K[x_1,x_2,\ldots,x_n]$ the polynomial ring over $K$ 
with $\deg(x_i)=1$ for all $i$. 
Let $e_i=e_i(x_1,\ldots,x_n)$ be 
the elementary symmetric function of degree $i$ in $R$, i.e., 
$$
e_i(x_1,\ldots,x_n)=\sum_{ j_1<j_2<\cdots<j_i}x_{j_1}x_{j_2}\cdots x_{j_i}
$$
for all $i=1,2,\ldots,n$. 

We denote by $S$ the subring  $S=K[e_1,e_2,\ldots,e_n]$ of $R$.  
Put $S_j=S \cap R_j$.  Then we have $S=\oplus_{j\geq 0}S_j$, which we regard as defining the grading of $S$, so 
that the natural injection $S \subset R$ is a grade-preserving algebra map. 
Let $\cH$ be the set of harmonic functions in $R$. 
Namely, $\cH$ is the vector space spanned by the partilal derivatives of the alternating polynomial 
$\prod _ {i < j}(x_i-x_j)$. 
The following are well known.  
\begin{itemize}
\item[(1)] 
The ring $S$ contains all symmetric functions of $R$.  
\item[(2)] 
$\cH$ is isomorphic to $R/(e_1, \ldots, e_n)$ as a graded vector space.  
\item[(3)]
The map 
$$
\cH\otimes_K S \ni h\otimes e \mapsto he \in R 
$$
is an isomorphism as graded vector spaces. 
Hence it follows that 
$R$ is a finite free module over $S$. 
\end{itemize}
\end{notation and remark}

\begin{lemma} \label{lemma:7-4}  %lemma7.5
With the same notation as above, 
let $f_1,f_2,\ldots,f_n\in S$ be homogeneous polynomials. 
Put $I=(f_1,f_2,\ldots,f_n)R$ and $J=(f_1,f_2,\ldots,f_n)S$. 
Then $R/I$ is a finite free module over $S/J$. 
\end{lemma}

\begin{proof}
This follows immediately from the last statement of Notation and Remark~\ref{NR:7-3}.  
\end{proof}

\begin{lemma} \label{lemma:7-5}   %lemma7.6 
The Artinian complete intersection 
$$
B=K[e_1,e_2,\ldots,e_n]/(p_a,p_{a+1},\ldots,p_{a+n-1})
$$ 
has the SLP, 
where $p_d=x_1^d+x_2^d+\cdots+x_n^d$ is the power sum of degree $d$. 
\end{lemma}

\begin{proof}
Since  
$$
h_B(q)=\frac{(1-q^a)(1-q^{a+1})\cdots(1-q^{a+n-1})}
{(1-q)(1-q^2)\cdots(1-q^n)} 
\ \ \ {\rm and}\ \ \  
q^{an-n}h_B(q^{-1})=h_B(q), 
$$
$B$ has a symmetric Hilbert function  
and 
the socle degree of $B$ is equal to $an-n$.  
Put 
$$
C=K[x_1,x_2,\ldots,x_n]/(x_1^a,x_2^a,\ldots,x_n^a). 
$$
One notices that $B$  and $C$  have the same socle degree. 

Next we show that $B$ is naturally a graded subring of   $C$.  
The symmetric group $G=S_n$ acts on $C$ 
by  permutation of the variables. 
Consider the exact sequence 
$$
0 \ra (x_1^a,\ldots,x_n^a) \ra R \ra C \ra 0.  
$$
Since ${\rm char}(K)=0$, 
we have the exact sequence
$$
0 \ra (x_1^a,\ldots,x_n^a)^G \ra R^G \ra C^G \ra 0, 
$$
where $M^G$  denotes the invariant subspace  for any $G$-module $M$.  
Note that $R^G=S$. 
Hence it follows that $C^G\cong S/((x_1^a,\ldots,x_n^a)\cap S)$. 
We would like to prove that 
$$(x_1^a,\ldots,x_n^a)\cap S=(p_a,p_{a+1},\ldots,p_{a+n-1}),$$ 
or equivalently, 
the natural surjection 
$$
\psi: B \ra  S/((x_1^a,\ldots,x_n^a)\cap S) \cong C^G 
$$
is an  isomorphims. 
By way of contradiction assume  that $\Ker(\psi)\neq(0)$. 
Since  $B$ is Gorenstein, $B_{an-n}$ is the unique  minimal ideal of $B$ and it  should be contained in 
$\Ker(\psi)$.    Hence the socle degree of $B/\Ker(\psi)$ is less than $an-n$. 
However, since the element $(x_1x_2\cdots x_n)^{a-1} \in C$  lies in  $C^G$, the socle degree of $C^G$ is 
equal to $an-n$. This is a contradiction.  
We have proved  $\Ker(\psi)=(0)$. 

It is known that the image of $e_1$ in $C$ is a strong Lefschetz element for $C$. (cf. Corollary 3.5 in \cite{jW87a}.)
Thus 
it follows that $e_1$ is a strong Lefschetz element for $B$   by Lemma~\ref{lemma:7-2}.  
\end{proof}

\begin{proof}
[Proof of Proposition~\ref{main-th4}]
The first assertion is proved in Lemma~\ref{lemma:7-1}.    
By Example 6.4 of \cite{tHjW06}, 
the algebra  
$$
A/mA=K[x_1,x_2,\ldots,x_n]/(e_1,e_2,\ldots,e_n)
$$ 
has the SLP.   
Hence the second  assertion follows 
from Lemmas~\ref{lemma:7-4},~\ref{lemma:7-5} 
and Theorem~\ref{main-th3}. 
\end{proof}

%%%%%%%%%%%%%%%%%%%%%%%%%%%%%%%%%%%%%%%%%%%%%%%%%%%%%%%%%%%%%%
%
%
%%%%      End of section 7 --- start of section 8     %%%%%%
%
%
%%%%%%%%%%%%%%%%%%%%%%%%%%%%%%%%%%%%%%%%%%%%%%%%%%%%%%%%%%%%%%

\section{Some more applications}   %section8  
Throughout this section we fix $R, S$ to be the same as in {\rm Notation and Remark~\ref{NR:7-3}}.  
Suppose that $f_1, \ldots, f_n$ is a  regular sequence of $S$.  
Put $B=S/(f_1, \ldots, f_n)S$ and $A=R/(f_1, \ldots, f_n)R$.
Let $m$ be the maximal ideal of $B$.  
Then we have (1) $A$ is finite flat over $B$ and (2) $A/mA$ has the SLP.  
Thus, by Theorem~\ref{main-th3}, if $B$ has the  SLP, then $A$ has the SLP.  
This is the idea of the proof of Proposition~\ref{main-th4}.  
We give  more applications of Theorem~\ref{main-th3}.

\begin{proposition} %prop8.1  
Let  $f \in S$ be  
a homogeneous element of degree $d$. 
Suppose that  
  $(e_2, e_3, \ldots, e_{n}, f)$ is a complete intersection in $S$.  
Then $R/(e_2, e_3, \ldots, e_{n}, f)R$ has the SLP.
\end{proposition}

\begin{proof}

Put $A=R/(e_2, \ldots, e_n, f)R$ and $B=S/(e_2, \ldots, e_n, f)$. 
It is easy to see that $A$ is finite flat over $B$.  
Note that $B$ has embedding dimension one. 
Hence it has the SLP, as one easily notices. 
Since the fiber of $B \ra A$ is  $R/(e_1, \ldots, e_n)R$, it has the SLP.
By Theorem~\ref{main-th3} the assertion is proved.

\end{proof}

\begin{lemma} \label{new_lemma_2}   %lemma8.2 
Suppose that  
$$
f_2, f_3, \ldots, f_n, f_d 
$$
is  a regular sequcence in $S$  such that $\mbox{\rm degree }f_i=i$ for 
$i=2, 3, \ldots, n, d$  with  $d > n$. 
Put $B=S/(f_2, \ldots, f_n, f_d)$.  
Then the following conditions are equivalent.   
\begin{enumerate}
\item 
$B$ has the strong Lefschetz property.
\item 
The embedding dimension of $B$ is 1. 

\end{enumerate}
\end{lemma}
\begin{proof}
Since  $B$  is a complete intersection, we have 
$$h_B(q)=1+q+ \cdots + q^{d-1}.$$ 
Assume (2).  Then as a $K$-algebra, $B$ is generated by $\ol{e_1}$.  
Hence $B$  is isomorphic to $K[X]/(X^d)$. Thus (1) follows. 
Assume (1).    The maximal ideal of $B$ is generated by 
$\ol{e_1}, \ol{e_2}, \cdots, \ol{e_n}$.  
Since $B$  has the strong Lefschetz property, we get that 
$\overline{e_k}$ is a constant multiple of $\overline{e_1}^k$.
Hence the maximal ideal is generated by a single element. 
\end{proof}

\begin{proposition}  %prop8.3 
Using   the same notation as the previous lemma, 
assume that  $d$  is a prime number.   
Then we have  
\begin{enumerate}
\item
$B=S/(f_2, f_3, \ldots, f_n, f_d)S$ has the strong Lefschetz property.  
\item
$A=R/(f_2, f_3, \ldots,  f_n,  f_d)R$ has the strong Lefschetz property.  
\end{enumerate}
\end{proposition}
\begin{proof}
(2) follows form (1)  by Theorem~\ref{main-th3}. 
To prove (1) it suffices to show that the embedding dimension of $B$ 
is one  by  Lemma~\ref{new_lemma_2}.  
Write $f_2= \al e_1^2 + \be e_2$, with $\al, \be \in K$.  
Assume $\be = 0$.  Then since 
$(e_1^2, f_3, \ldots, f_n, f_{d})$ is a complete intersection, 
so is $(e_1, f_3, \ldots, f_n, f_{d})$.
Hence, if we put $B'= S/(e_1, f_3, \ldots, f_n, f_{d})$, then 
\[
h_{B'}(q)=\frac{(1-q)(1-q^3)\cdots (1-q^n)(1-q^d)}{(1-q)(1-q^2)\cdots (1-q^n)}
=\frac{1+q+ \cdots +q^{d-1}}{1-q}.
\]
This forces the numerator of the last function be  divisible by $1-q$, a contradiction, 
since $d$ is a prime. 
This means that $\be \neq 0$.  
Thus we may replace  $e_2$ by  $f_2$ as a generator of the 
algebra $B$.   Hence, by modifying the elements $f_3, \ldots, f_n, f_d$ suitably, 
we have   
\[
B \cong K[e_1, e_3, \ldots, e_n]/(f_3, \ldots, f_n, f_d).  
\]
Now suppose that $f_3= \al e_1^3 + \be e_3$.  As with the preceeding case it does not 
occur that $\be = 0$.  Hence $f_3$ may replace $e_3$  as a generator of the algebra.  
Hence we have 
\[
B \cong K[e_1, e_4, \ldots, e_n]/(f_4, \ldots, f_n, f_d), 
\]  
with modification of the generators  $f_4, \ldots, f_n, f_d$. 
We may repeat the same argument to obtain 
\[
B \cong K[e_1]/(e_1 ^d). 
\]
\end{proof}
%%%%%%%%%%%%%%%%%%%%

\begin{remark}   %rem8.4  
Denote by $[a]_q$ the 
$q$-integer $[a]_q=1+q+ \cdots +q^{a-1}$.  
Then in the statement of the above proposition, it is the same if we say 
 $[d]_q$ is not divisible by   any one of  $[2]_q, \ldots, [n]_q$, 
instead of 
the assumption ``$d$ is a prime number.'' 
\end{remark}

%%%%%%%%%%%%%%%%%%%%

\section{Appendix}

We show a new proof for the following 

\begin{proposition}  \label{ikeda}   %%%prop9.1 

$K[X,Y]/(X^r, Y^s)$ has the strong Lefschetz property with $X+Y$ a strong Lefschetz element, 
where $\mbox{\rm char}\  K=0$ or $\geq {\rm Max}\{r,s\}$ and $\deg \ X = \deg \ Y =1$.  

\end{proposition}

\begin{proof}   
H. Ikeda~\cite{hI96} proved the following: 
If $A$ is a standard graded  ring and if  $x \in A$ is a strong Lefschetz element for $A$, then 
$A[y]/(y^2)$ has the strong Lefschetz property with 
$x+\ol{y}$ a strong  Lefschetz element.  (In fact, this is a special case of Theorem~2.8 in \cite{hI96}.)
We would like to emphasize  that her proof does not use the theory of $sl_2$ or the theory of  Groebner bases. 
Using her thoerem  we  immediately get  that 
$$A:=K[z_1, \ldots, z_n]/(z_1^2, \ldots, z_n^2)$$
has the strong Lefschetz property with $\ol{z}_1 + \cdots + \ol{z}_n$  a  
strong Lefschetz element.  
Now choose $n$ so that $n=r+s-2$ and consider the ring homomorphism
\[
\phi : K[X,Y] \ra  A 
\]
defined by $X \mapsto \ol{z}_1+ \cdots +  \ol{z}_{r-1}$ and $Y \mapsto \ol{z}_{r} + \cdots + \ol{z}_{n}$.  
Put $J= \ker \  \phi$.  
We claim that $J=(X^{r}, Y^{s})$.   
It is easy to see the inclusion 
\begin{equation}  \label{natural_map} 
J \supset (X^r, Y^s). 
\end{equation}
We have to show that it is the equality.  By way of contradiction assume that $J$ is strictly larger than 
$(X^r, Y^s)$.  
Then the homogeneous part of degree $r+s-2$  of $R/J$ is $0$, since $R/(X^r, Y^s)$ has the 
unique minimal ideal, the socle,  at that degree.  
On the other hand we have  
$$\phi(X^{r-1}Y^{s-1})  = (r-1)!(s-1)!\; \ol{z_1z_2 \cdots z_n} \neq  0.$$
This is a contradiction.      
Thus  the inclusion (\ref{natural_map})  is in fact the equality. 
We have just shown  that the subring  $B:=K[\phi(X),\phi(Y)]$ of $A$  is isomorphic to 
$$K[X,Y]/(X^r, Y^s).$$ 
Note that $B$ has a symmetric Hilbert function and moreover  
shares  the same socle with $A$.    Then,  
since the Lefschetz element 
$\ol{z}_1+ \cdots + \ol{z}_{n}$ for $A$ lies in  $B$, it is a Lefschetz element for $B$ also.   
\end{proof}

\end{document}